\documentclass[a4paper,10pt]{article}
\usepackage{geometry,times}
\usepackage{amsfonts}
\usepackage{amssymb}
\usepackage{amsmath}
\usepackage{tikz}
\usepackage{pgfplots}
\usepackage{tikz-3dplot}
\tdplotsetmaincoords{60}{115}
\pgfplotsset{compat=newest}
\definecolor{paleyellow}{rgb}{1, 1, 0.9}
\usepackage{amsthm}
\newcommand{\rose}{R\accent 23 u\v zi\v cka}
\newcommand{\bb}{\mathbb}

\theoremstyle{plain}
\newtheorem{satz}{prop}
\newtheorem{corollary}[satz]{Corollary}

\newtheorem{definition}[satz]{Definition}

\theoremstyle{thm}
\newtheorem{theorem}[satz]{Theorem}
\theoremstyle{note}

\theoremstyle{problem}

\theoremstyle{notadef}

\usepackage[]{graphicx}
\usepackage{color}
\definecolor{dark-BLGR}{rgb}{0.28,0.28,0.28}

\definecolor{dark-red}{rgb}{0.90,0.14,0.14}

\definecolor{dark-viol}{rgb}{0.73,0.18,0.69}
\definecolor{Cyan}{rgb}{0.31,0.67,0.82}
\newcommand{\blue}{\textcolor{blue}}
\newcommand{\red}{\textcolor{red}}
\newcommand{\green}{\textcolor{green}}

\definecolor{light-yellow}{rgb}{1,1,0.8}

\definecolor{Blue}{rgb}{0.0353,0.0275,0.4}

\definecolor{GRAY}{rgb}{0.26,0.26,0.26}

\definecolor{GREYY}{rgb}{0.45,0.45,0.45}

\definecolor{vivid-viol}{rgb}{0.3255,0.0353,0.55}

\definecolor{dark-blue}{rgb}{0.05,0.05,0.65}

\definecolor{dark-green}{rgb}{0.03,0.77,0.29}

\definecolor{dark-Green}{rgb}{0.03,0.57,0.09}

\definecolor{strong-viol}{rgb}{0.2353,0.094,0.349}

\definecolor{BLUE}{rgb}{0.41,0.44,0.93}

\definecolor{RED}{rgb}{0.90,0.18,0.32}

\definecolor{Yel}{rgb}{0.89,0.95,0.19}

\definecolor{White}{rgb}{1,1,1}
\definecolor{black}{rgb}{0,0,0}

\definecolor{Orchid}{rgb}{0.6,0.1607,0.2}

\definecolor{Orange}{rgb}{0.99,0.49,0}

\definecolor{Thistle}{rgb}{0.1151,0.1249,0.9873}

\definecolor{brown}{rgb}{0.298,0.153,0.07843}

 \usepackage{titlesec}
\titleformat{\section}
{\normalfont
  \bfseries
}
{\thesection}{1ex}{}
\titleformat{\subsection}
{\normalfont
  \bfseries
}
    {\thesubsection}{1ex}{}

\newcommand{\cA}{{\mathcal{A}}}

\begin{document}
\hspace*{1cm}\\[.4cm]
{\large \textbf{Exact Poincaré Constants in
  three-dimensional Annuli}}
\\[.5cm]
Bernd Rummler\\
\hspace*{2mm}{\small \em Otto-von-Guericke-Universit\"at Magdeburg,
  Inst. f\"ur Analysis und Numerik, PF
  4120, 39016 Magdeburg}
\\[2mm]
Michael \rose\\
\hspace*{2mm}{\small \em Inst. f{\"u}r Angewandte Mathematik, Universit{\"a}t Freiburg,
  Ernst-Zermelo-Str.~1, 79104 Freiburg}
\\[2mm]
Gudrun Th\"ater\footnote{Corresponding author\quad
  E-mail:~\textsf{gudrun.thaeter@kit.edu}}\\
\hspace*{2mm}{\small \em Inst. f\"ur Angewandte \& Numerische Mathematik, KIT, 
  76128 Karlsruhe}
\\[5mm]
\hspace*{.1cm}\hfill\parbox{16cm}{
{\small\textbf{Abstract:}
 We study 3d-annuli. In our non-dimensional setting each annulus ${\Omega}_{\cal A}$ is
defined via two concentrical 
balls with radii ${\cal A}/2$ and ${\cal A}/2 +1$. For these
geometries we provide the exact value for the Poincar\'e constants for
scalar functions and calculate precise Poincar\'e constants
for solenoidal 
vector fields 
(in both cases with vanishing Dirichlet traces on the
boundary). For this we use the 
first eigenvalues of the scalar Laplacian and
the Stokes operator, respectively. 
Additionally, corresponding problems in domains
${\Omega}_{\sigma}^{*}$, the 3d-annuli from \cite{RumTh2024},
are investigated - for comparison but also to provide limits for  ${\cal
  A}\,\to\,0$.  
In particular, the Green's function of the Laplacian on
${\Omega}_{\sigma}^{*}$ with vanishing 
Dirichlet traces on $\partial {\Omega}_{\sigma}^{*}$ is used
to show that for ${\sigma}\,\to\,0$ the first eigenvalue here 
tends to the first eigenvalue of the corresponding problem on the open
unit ball. 
On the other hand, we take advantage of the so-called small-gap limit
for  ${\cal A}\to\infty$. 
}}
\newcommand{\abs}[1]{\left | #1 \right |}
\newcommand{\RR}{\mathbb{R}}
\newcommand{\cC}{{\mathcal{C}}}
\section{Introduction\label{sec_int}}
In modelling the behaviour of very diverse objects such as planets or thin layers
around balls one has to understand 
fluid flow in spherical shells - i.e. flow in the domain between
two concentric spheres (the configuration  
is depicted in Fig.~\ref{fig:geometry} below).
\begin{figure}[th]
\centering
\begin{tikzpicture}[scale=2,tdplot_main_coords] 
\coordinate (P) at ({0},{1.27/sqrt(2)},{1.27/sqrt(2)
});
\coordinate (P1) at ({0},{-0.51/sqrt(2)},{0.51/sqrt(2)});
\shade[ball color = blue,
    opacity = 0.25
] (0,0,0) circle (1cm);
\shade[ball color = gray,
    opacity = 0.3
] (0,0,0) circle (.5cm);
 
 
 
\draw[-stealth] (0,0,0) -- (1.80,0,0);
\draw[-stealth] (0,0,0) -- (0,1.30,0);
\draw[-stealth] (0,0,0) -- (0,0,1.30);
 
\draw[thick,blue, -stealth] (0,0,0) -- (P);
\draw(0.5,1.0,1.2) node[anchor=east,blue]{$R_o$};
\draw[thick,gray, -stealth] (0,0,0) -- (P1);
\draw((0.25,-0.05,0.35) node[anchor=east]{$R_i$};
%
\end{tikzpicture}
\caption{3d-annulus  with inner/outer radius  in a
Cartesian coordinate system}
\label{fig:geometry}
\end{figure}
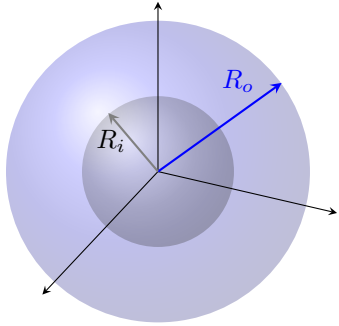
Here one often has to model {\em Natural convection}, i.e.~the
fluid moves (i) under the influence of gravity and (ii) because of
temperature differences causing changes in the density of the
fluid. This is typical, e.g., for astrophysical resp. geophysical
applications. 
From the geophysical
point of view the flow within a spherical gap provides a model for
movement in the atmosphere or within the Earth’s mantle.  
The question of the origin of the Earth’s magnetic field -- generated
by currents in the outer core -- has recently been explored with
increasingly realistic parameters, enabled by advances in
computational power. 
Additionally, this system is of interest from a fundamental
perspective on pattern formation. Many fluid-mechanical systems
exhibit symmetry-breaking flows due to various hydrodynamic
instabilities, which often produce intriguing geometric
shapes. These systems have led to significant insights and the
development of broadly applicable analytical methods (see, e.g. 
\cite{Junk}). 
The possible flow patterns are very complex and allow, depending on
the material properties 
(encoded in parameters as the Rayleigh and Prandtl number, e.g.) 
as well as the actual relation of the radii,
very diverse behaviour. It is interesting to find
stability regions for each of the observed patterns.
Note, that precise knowledge about Poincar\'e constants 
(see Eq.~\eqref{Poinc}  below) makes all
results for, e.g., existence, uniqueness and stability of solutions
more precise. 

In our three-dimensional setting here we will apply  
similar methods as we developed for two-dimensional
annuli. 
They can be studied as
  the cross-section of the domain between  
horizontal coaxial cylinders and from experiments one knows that for a 
broad range of parameters  in this
configuration the flow behaviour is mainly
two-dimensional.\footnote{Investigations of convective flow between
  horizontal coaxial cylinders in 
a gravitational field show that 
for small Rayleigh numbers and for any relation of the radii there exists a
steady 2d-flow with two crescent-shaped eddies which are symmetric with
respect to the vertical plane through the common axis of the
cylinders. 
As the Rayleigh number increases above a critical value, 
different flow patterns are observed. To understand the stability
behaviour and
to find these bifurcations are very delicate tasks. 
}
Thus, the study of the cross-section provides 
the information for the three-dimensional flow problem in that
case.
In particular, in \cite{RuRuTh2016} 
we provide optimal analytical bounds for the Poincar\'e
constants 
as functions of the non-dimensional parameter ${\cal A}$ (defined in
(1) below).

An approach to narrow down possible regions for bifurcations
can be found in \cite{passerini2009,passerini2010} and references
therein. It is very interesting for the  study of
stability to
find the critical Rayleigh number as a function of the radii
relation in an analytical way. A very encouraging result in this
direction is described in \cite{passerini2024}.
\\[3mm] 
In the present paper we calculate the Poincar\'e
constants
as function of 
a non-dimensional number, which characterises the (relative) size of
the gap.
There are two frames:  Using the so-called {\em
  inverse relative gap width} (mathematically more precise would be 
{\em gap weighted inner diameter}) $\cA$ or the parameter ${\sigma}$. More
precisely, for $R_{i}$ and $R_{o}$  denoting the inner and
outer radius of the annulus, respectively, and $0\,<\,R_{i}<\,R_{o}$
(see Fig.~\ref{fig:geometry}) both quantities are defined as
\begin{equation}\label{calA}
{\cal A}\,:=\,\frac{\displaystyle{2 R_{i}}}{\displaystyle{ R_{o}
    -R_{i}}}
\qquad\mbox{while}\qquad
{\sigma}\,:=\,\frac{\displaystyle{R_{i}}}{\displaystyle{ R_{o}}}=
\frac{\displaystyle{\cal A}}{\displaystyle{ {\cal A} +2 }}\,.
\end{equation}
Our \textbf{main results} are: To provide the exact value for the Poincar\'e constants
for scalar fields as
well as very good approximations for 
solenoidal vector fields 
{\em including the limits} for the extreme cases, i.e. infinitely large and zero gap width.
\\[3mm]
{\bf General notation A.} Let ${\bb R}^{3}$ be endowed with the
usual Euclidian norm $\| . \|$ 
and  elements of ${\bb R}^{3}$ be denoted by underlined small letters.
The open unit ball is ${\Omega}^{*}:=\{{\underline{x}} \in {\bb R}^{3}:\,
\|{\underline{x}}\|\,
<1\}$, the unit spherical surface  $\omega\,:=\,\{{\underline{x}}\in {\bb
  R}^{3}:\,\|{\underline{x}}\|=1\}$ 
and closed spherical surfaces around the origin with radius $r$ are
$\omega_{r}\,:=\,\{{\underline{x}}\in {\bb
  R}^{3}:\,\|{\underline{x}}\|=r\} 
=r\omega$
for all $r\,\in\,(0,\infty)$.
\\[3mm]
{\bf Annulus domains.} It is useful to study our domain without
dimensions. As usual 
we pick the annulus (the sherical shell) with fixed gap width $1$ using our non-dimensional
parameter 
${\cal A}$ as follows:
For any ${\cal A}\,\in\,(0,\infty)$, we denote
by
\[
{\Omega}_{{\cal A}}\,:=\{{\underline{x}} \in {\bb R}^{3}:\,{\cal A}/{2}
<\|{\underline{x}}\|<1+{\cal A}/{2}\}\,.
\]
Its boundary  
$\partial{\Omega}_{\cal A}$ consists of the two parts
$\omega_{{\cal A}/{2}}$ and  $\omega_{1+{\cal A}/{2}}\,.$  Moreover,
 for all   
${\sigma}\,\in\,(0,1)$ we introduce the family of annuli \[
{\Omega}^{*}_{\sigma}:=\{{\underline{x}} \in {\bb
  R}^{3}:\,0<\sigma<\|{\underline{x}}\|\,<1\}\,.
\]
With this notation we can directly use results from \cite{RumTh2024} and 
it is more convenient for the case
${\cal A}\,\to\,0$.
The advantage is that these three-dimensional annuli are subsets of the
unit ball 
${\Omega}^{*}$. They have  the boundary
$\partial{\Omega}^{*}_{\sigma}=
\omega_{\sigma} \cup \omega$. 
\\[3mm]
{\bf General notation B.} Let ${\Omega}$ stand as shorthand for any of the domains
${\Omega}_{{\cal A}},\, {\Omega}^{*}_{\sigma},$ and ${\Omega}^{*}$
and the abbreviation $(.)$ for $({\Omega})$, respectively. 
We consider the usual Lebesgue and Sobolev spaces ${\bb L}_{2}(.)$ and
${\bb W}_{2}^{k}(.)$ 
of scalar functions and 
${\underline{\bb L}}_{{\;\!}{2}}(.)=({\bb L}_{2}(.))^{3}$
and ${\underline{\bb W}}_{{\;\!}{2}}^{k}(.)=({\bb
  W}_{2}^{k}(.))^{3}$ of vector functions. 
The norm in ${\bb L}_{2}(.)$ is denoted by $\| . \|_{2}$,
${\bb W}_{2}^{1}\hspace{-.62cm}{~}^{{~}^{{~}^{o}}}\hspace{.2cm}(.)$ is
the closure of $C_{o}^{\infty}(.)$ 
in ${\bb W}_{2}^{1}(.)$.  
All solenoidal vector functions belonging to 
${\underline{C}}_{{\;\!}{o}}^{\infty}(.)$ form $\underline{\cal V}(.)$. The closures
of $\underline{\cal V}(.)$ in ${\underline{\bb L}}_{{\;\!}{2}}(.)$ 
and
${\underline{\bb W}}_{{\;\!}{2}}^{1}(.)$, respectively, are denoted by
${\underline{\bb H}}(.)$ 
and ${\underline{\bb V}}(.)$, respectively.  We use the spherical
Bessel functions ${J_{k}(.)}$ of the first kind as well as 
the spherical Bessel functions  ${J_{-k}(.)}$ 
of order  $k\,\in\, \{\frac{1}{2}+m,\,m\,\in\, {\bb N}_{o}\} $. The 
polar coordinates are $r$,  $\vartheta$ and $\varphi$ with the corresponding unit vectors
${\underline{\mathfrak{e}}}_{{\;\!}r}$, ${\underline{\mathfrak{e}}}_{{\;\!}\vartheta}$ and
${\underline{\mathfrak{e}}}_{{\;\!}\varphi}$.
\\[3mm]
The {\bf Poincar\'e-(Friedrichs-)inequalities} are the central tools
to ensure, that  the spaces ${\bb
  W}_{2}^{1}\hspace{-.62cm}{~}^{{~}^{{~}^{o}}}\hspace{.2cm}(.)$ and  
${\underline{\bb V}}(.)$ 
can be equipped with equivalent
norms generated by the Dirichlet norms 
(which otherwise would only be semi-norms):
\begin{align} \label{Dirichlet} 
\|u\|_{D}  :=
\Big(\sum_{k=1}^{3}\Big\|{\frac{\displaystyle{\partial
  u}}{\displaystyle{\partial x_{k}}}}\Big\|_{2}^{2} 
\Big)^{1/2}\quad \forall \, u\,\in\,{\bb
  W}_{2}^{1}\hspace{-.62cm}{~}^{{~}^{{~}^{o}}}\hspace{.2cm}(.)\,,\quad 
\|{\underline{u}}\|_{D,S}   := 
\Big(\sum_{j,k=1}^{3}
\Big\|
{\frac{\displaystyle{\partial u_{j}}}{\displaystyle{\partial x_{k}}}}
\Big\|_{2}^{2}
\Big)^{1/2}
\quad \forall \,  {\underline{u}}  \,\in\,{\underline{\bb V}}(.)\,,
\end{align}
where the so-called Frobenius inner product is involved in the last
definition. 
Denoting by $c_{p}({\cal A})$ and $c_{p,S}({\cal A})$ the Poincar\'e
constants with respect to the spaces ${\bb
  W}_{2}^{1}\hspace{-.62cm}{~}^{{~}^{{~}^{o}}}\hspace{.2cm}({\Omega}_{{\cal
    A}})$ and ${\underline{\bb V}}({\Omega}_{{\cal A}})$,
respectively, the Poincar\'e-(Friedrichs-)inequalities are 
\begin{align}\label{Poinc}
\|u\|_{2}  \leq   c_{p}({\cal A})   \|u\|_{D} \quad\forall \, u\,
\in\,{\bb W}_{2}^{1}\hspace{-.62cm}{~}^{{~}^{{~}^{o}}}
\hspace{.2cm}({\Omega}_{\cal A})\,\quad
\mbox{and}\quad
\|{\underline{u}}\|_{{\underline{\bb L}}_{{\;\!}{2}}}
\leq c_{p,S}({\cal A})
\|{\underline{u}}\|_{D,S} \quad
\forall \,  {\underline{u}}  \,\in\,
{\underline{\bb V}}({\Omega}_{\cal A})\,.
\end{align}
The Poincar\'e constants 
are related to the 
first eigenvalue of the Laplace or Stokes operator on the
${\Omega}_{\cal A}$-domains 
(with vanishing Dirichlet traces), respectively,
because of the relations (see the Ths. in Subsections 4.5.3 and 4.5.4 and Th. 3 in
6.1.5 in \cite{Triebel})
\begin{equation} \label{simple-tool}
c_{p}({\cal A}) = 
\big(\lambda_{1,L}({\cal A})\big)^{-1/2}
\quad
\,\quad \mbox{and}\quad
c_{p,S}({\cal A}) =
\big(\lambda_{1,S}({\cal A}) \big)^{-1/2} \,,
\end{equation}
where ${\lambda}_{1,L}({\cal A})$ and  ${\lambda}_{1,S}({\cal A})$
denote the first simple eigenvalue of the Laplace operator and the
first (triple) eigenvalue of the Stokes  
operator, respectively. One important result in our setting is that
$\forall \,{\cal A} \,\in\, [0,\infty)$ we find ${\lambda}_{1,L}({\cal A})\,=\,\pi^2$. 
We refer to \cite{RumTh2024} for very detailed 
information with respect to our domain and  to \cite[Subsection 6.4.4]{Triebel} for the
open unit ball ${\Omega}^{*}$. 
\\[.3cm]
Our paper is organised as follows:
We collect the essential theoretical fundamentals in
Section \ref{Sec2}. There we sketch 
the procedures to construct the Laplace as well as the Stokes operator
as Friedrichs' extension from the Poisson and the Stokes
problem, respectively. 
We benefit from the properties of operators with a pure real
point spectrum. 
We introduce the {\em Leray-Helmholtz projector} 
$\Upsilon : {\underline{\bb
    L}}_{{\;\!}2}(.)\,\longmapsto \,{\underline{\bb H}}(.)$ 
as well as the criteria for the smallest eigenvalues
${\lambda}_{1,L}({\cal A})$, ${\lambda}_{1,S}({\cal A})$, 
and their corresponding eigenfunctions.
Section \ref{Inv_Lim}  is devoted to the limiting cases. We carefully
conduct the transition ${\cal A}\,\to\,0$ especially for 
${\lambda}_{1,L}({\cal A}\,\to\,0)$ in the form of
${\lambda}_{1,L}({\sigma}\,\to\,0)$. 
The investigations are performed with the Green's functions for
circular annuli ${\Omega}^{*}_{\sigma}$ 
and for the unit ball ${\Omega}^{*}$.  The crucial result here is,
that as $\sigma\,\to\,0$ the  {\it{problem forgets}} the 
center point together with the boundary condition there.

The study of the behaviour for ${\lambda}_{1,S}({\cal A}\,\to\,0)$ is
much easier, because of the  
vanishing of the first eigenfunction of the Stokes operator for
${\Omega}^{*}$ in  
${\underline{x}}={\underline{0}}$.
Finally, by simple transformations we show, that the cases
${\lambda}_{1,L}({\cal A}\,\to\,\infty)$ 
and ${\lambda}_{1,S}({\cal A}\,\to\,\infty)$ are covered by the
so-called {\em small gap limit.} In Section \ref{Sec4} the values of
the Poincar\'e constants   
$c_{p}({\cal A})\,=\,{\frac{1}{\pi}}$ 
and the calculated values for the Poincar\'e constants $c_{p,S}({\cal A})$ are 
represented as a graph  for
${\cal A}\in [0,\infty)$. 
\section{Theoretical groundwork}\label{Sec2}
\subsection{Available Bounds for the Poincaré constant}\label{Sec2avBou}
As rule of thumb it holds: The Poincaré constant 
depends on the diameter
of the considered domain $\Omega$ (cf. \cite{RuRuTh2016}). A familiar estimate 
for $\Omega\subset \bb R^3$ is $  c_p \leq {\text{diam}(\Omega)}/2$
(see, e.g., \cite[Th.~4.1]{galdi1998}).

A certain improvement can be achieved by
$ c_p \leq {\text{diam}({\Omega}_{\cal A})}/( {\pi\sqrt 2})$ in Exercise~4.2 therein. For
$\Omega_{\mathcal A}$ and with the abbreviation $\kappa_{1,L}({\cal
  A})$ for the smallest positive solutions 
of the transcendental equation \eqref{trans1} (respectively $\kappa_{1,L}({0})$ for the smallest positive root of 
$J_{\frac{1}{2}}(.)$ in sense of the limit ${\cal A}\,\to\,0 $) one finds:
\begin{equation*}
c_p \,=\,{\frac{1}{\kappa_{1,L}({.})}}
\leq \frac{\sqrt 2}{\pi}
\left( 1+\frac{\mathcal A}{2}\right)
\,.
\end{equation*}
By using a result from \cite{nazarov2000}, we obtain the following
(improved) bound
\begin{equation*}
\label{eq:bound_new}
  c_p \leq {\frac{R_o}{R_i}} \, \frac{1}{\pi}
=\frac{1}{\pi} \, \left( {1 + \frac{2}{\mathcal A}}\right)\,.
\end{equation*}
\begin{proof}
With spherical polar coordinates
$\{\underline{r}:=(r,\vartheta,\varphi)^{T}\,:\,0\leq
r<\infty\,,\,0\leq\vartheta\leq\pi\,,\,0\leq\varphi < 2\pi\}$ 
for $\tilde w(\underline{r} 
) =
w({\underline{x}})$  (both
varying in their respective space ${\bb
  W}_{2}^{1}\hspace{-.62cm}{~}^{{~}^{{~}^{o}}}\hspace{.2cm}({\Omega}_{\cal
  A})$ ) we obtain that
\begin{align}
\label{eq:poincarepolar}
\hspace*{-1cm}
   c_p^2 =\max_w \frac{\int_{{\Omega}_{\cal A}} \left|w\right|^2 \,
d\underline{x}}{\int_{{\Omega}_{\cal A}} (\underline{\nabla} w)^T \cdot \underline{\nabla} w \, d\underline{x}}
= \max_{\tilde{w}} \frac{\int_{{\Omega}_{\cal A}} r^2 \sin
  (\vartheta)\left|\tilde w\right|^2 \, d\vartheta d\varphi dr
  }
     {\int_{{\Omega}_{\cal A}} \!\! r^2 \sin (\vartheta) \left( \left(\partial_r \tilde w\right)^2 +
           r^{-2} \left( \partial_\vartheta \tilde w \right)^2 +
           r^{-2} \sin^{-2} (\vartheta)\left( \partial_\varphi \tilde w \right)^2\right)\,
  d\vartheta d\varphi dr
  }\,. 
\end{align}
We see from \cite[Prop.~1.1]{nazarov2000} that the argument for which
\eqref{eq:poincarepolar} attains its maximal value is a radial function
$\tilde u_1 = \tilde u_1(r)$. Hence,  
\begin{equation}
\label{eq:poincarepolarradial}
   c_p^2 = {\int_{\frac{\mathcal A}{2}}^{1+\frac{\mathcal A}{2}} r^2 \abs{\tilde u_1(r)}^2 \,
     dr}\Big/{\int_{\frac{\mathcal A}{2}}^{1+\frac{\mathcal A}{2}} r^2 \abs{\tilde u_1'(r)}^2 \, dr}. 
\end{equation}
Note that $\tilde u_1$ is the eigenfunction associated to the smallest
eigenvalue of the Laplace problem in $\Omega_A$. 
Further, for the one-dimensional
Laplace eigenvalue problem we know from \cite{nazarov2000} 
that if 
$\tilde w$ varies  in $ {\bb
  W}_{2}^{1}\hspace{-.62cm}{~}^{{~}^{{~}^{o}}}\hspace{.2cm}(\frac{\mathcal
  A}{2},1+\frac{\mathcal A}{2})$ we may conclude
\begin{equation}
\label{eq:poincareEins}
\left(\frac{
\frac{\mathcal A}{2}+1-\frac{\mathcal A}{2}}{\pi}\right)^2 
=
\left(\frac{1}{\pi}\right)^2 = 
\max_{\tilde w } \frac{
\int_{\frac{\mathcal A}{2}}^{1+\frac{\mathcal A}{2}}
\abs{\tilde w(r)}^2 \, dr}
{
\int_{\frac{\mathcal A}{2}}^{1+\frac{\mathcal A}{2}}
\abs{\tilde w'(r)}^2 \, dr}
\end{equation}
\begin{equation}\label{eq:poincareZwei}
\,\Rightarrow\,
c_p^2 \leq \left(
\frac{1+\frac{\mathcal A}{2}}{\frac{\mathcal A}{2}}\right)^2
 \, \frac{\int_{\frac{\mathcal A}{2}}^{1+\frac{\mathcal A}{2}}
 \abs{\tilde u_1(r)}^2 \, dr}
 {\int_{\frac{\mathcal A}{2}}^{1+\frac{\mathcal A}{2}} \abs{\tilde u_1'(r)}^2 \, dr}         \leq  \left(
\frac{1+\frac{\mathcal A}{2}}{\frac{\mathcal A}{2}}\right)^2\, \left(\frac{1}{\pi}\right)^2\,,
\end{equation}
which is implied by \eqref{eq:poincarepolarradial}.
\end{proof}
\noindent As we see, this bound is useful (i.e.~small) for large
$\mathcal A$. 
So the up to now best available result is
\begin{equation}\label{bestpoincare}
  c_p \leq \min \Big\{ \frac{\sqrt 2}{\pi} \Big( 1+\frac{\mathcal
    A}{2} \Big) \,,
\frac{1}{\pi} ({1+\frac{2}{\mathcal A}}) \Big\}.
\end{equation}
\begin{figure}
\centering
\includegraphics[width=0.65\textwidth]{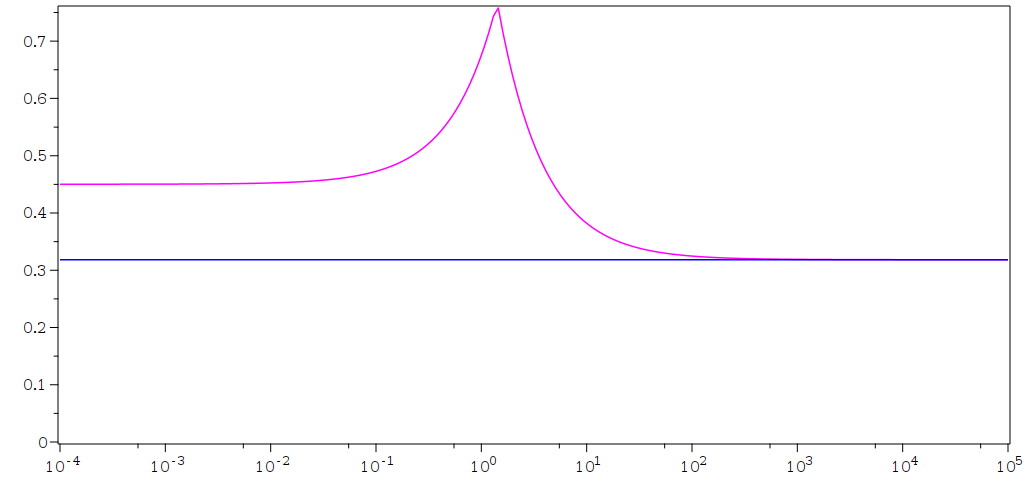}
\caption{Analytical bounds of Poincar\'{e} constants and  
  the constant value $\frac{1}{\pi}$ vs. $\mathcal A$}
\label{fig:poincare_plot}
\end{figure}
\hspace*{1cm}\\[-4mm]
The good news of Fig.~\ref{fig:poincare_plot} is that the fit of our
best analytical bound is fine for large $\cal A$. Unfortunately, for
small ${\cal A}<10$ the  
estimates are right but still poor. The constant $\frac{1}{\pi}$ is
very meaningful. 
Considering that for all ${\cal A}>0$ we know that
${\lambda}_{1,L}({\cal A})\,\leq\,{\lambda}_{1,S}({\cal A})$ since the
eigenvalue problem of the Stokes operator is a restricted one for the
Laplace operator. We see that $c_{p,S}({\cal A})\,\leq \,c_{p}({\cal
  A}) $, 
cf. \eqref{simple-tool},  and that \eqref{bestpoincare} is an estimate for the 
constant $c_{p,S}$ as well.
\subsection{Laplace and Stokes operators on
  3d-annuli} \label{sec_theo} 
As in General notation B, we take both symbols ${\Omega}$ and $(\cdot)$ as
placeholders.
\begin{definition}\label{D4} {\em The Laplace operator is
  defined as 
\begin{align*}
{\boldsymbol L^{\circledast}}\,{v} := -\Big(
{\displaystyle{\frac{{\partial}^{2} v}{\partial x_{1}^{2}}}}+{\displaystyle{\frac{{\partial}^{2} v}{\partial x_{2}^{2}}}}+
{\displaystyle{\frac{{\partial}^{2} v}{\partial x_{3}^{2}}}}\Big)
=- \Delta{\;\!}_{\underline{x} } v
\hspace*{0.7cm}  \forall\,v\,\in\,D({\boldsymbol
  L^{\circledast}})=C_{o}^{\infty}({\Omega})\,. 
\end{align*}
We denote the Friedrichs' extension of ${\boldsymbol L^{\circledast}}$
by ${\boldsymbol 
  L}:={\overline{\boldsymbol L^{\circledast}}}$, where ${\boldsymbol
  L}$ is applied on 
$D({\boldsymbol L})\,:=\,{\bb
  W}_{2}^{1}\hspace{-.62cm}{~}^{{~}^{{~}^{o}}}\hspace{.2cm}({\Omega}) 
\cup {\bb W}_{2}^{2}({\Omega})$.}
\end{definition}
\noindent{\bf Remark:} The range of ${\boldsymbol L}$ is
$R({\boldsymbol L})={\bb L}_{2}({\Omega})$. In this sense we may write:
${\boldsymbol L}=-\Delta{\;\!}_{\underline{x} } : D({\boldsymbol
  L})\,\longmapsto \,{\bb L}_{2}(.)$.
\\[.3cm] 
We need the Leray-Helmholtz projection $\Upsilon$ to define the Stokes operator. 
$\Upsilon$ is the well-defined
projector of ${\underline{\bb L}}_{{\;\!}2}(.)$ onto its subspace 
${\underline{\bb H}}(.)$
of generalised solenoidal fields with vanishing generalised traces in
the normal direction. 
We note, that 
it is also used in the sense of:
$\Upsilon :{\underline{\bb W}}_{{\;\!}2}^{1}(.)\,\longmapsto
\,{\underline{\bb V}}(.)$\,.
\begin{definition}\label{D5} 
{\em The Stokes operator is defined as
$
{\boldsymbol S^{\circledast}}\,{\underline{v}}:=-
\Delta_{\underline{x} } {\underline{v}}\hspace*{0.5cm}
\forall\,{\underline{v}}\in D({\boldsymbol
  S^{\circledast}})=\underline{\cal V}({\Omega})\,. 
$
We denote the Friedrichs' extension of ${\boldsymbol S^{\circledast}}$ by 
${\boldsymbol S}:={\overline{\boldsymbol S^{\circledast}}}$, where
${\boldsymbol S}$ is defined on its domain 
$D({\boldsymbol S}):={\underline{\bb S}}_{{\;\!}}^{2}(.)=
{\underline{\bb W}}_{{\;\!}2}^{2}(.)\cap{\underline{\bb V}}(.)$ \,.} 
\end{definition}
\noindent {\bf Remark:} The range of ${\boldsymbol S}$ is $R({\boldsymbol
  S})={\underline{\bb H}}(.) $. In this context one may write
${\boldsymbol S} =-\Upsilon \Delta{\;\!}_{\underline{x}
}:{\underline{\bb S}}_{{\;\!}}^{2}(.) 
\,\longmapsto \,{\underline{\bb H}}(.)$.
\\[.3cm]
We sketch the fundamental properties of both operators (i.e. ${\boldsymbol L}$
as well as ${\boldsymbol S}$) using ${\boldsymbol S}$ as example.
\begin{theorem} \label{thmstok}
The Stokes operator ${\boldsymbol S}$ is positive and self-adjoint.
Its inverse ${\boldsymbol S}^{-1}$ is injective, self-adjoint and compact.
\end{theorem}
\noindent The proof of Theorem \ref{thmstok} is a simple modification of 
Theorems 4.3 and 4.4 in \cite{CoFoi}. The essential tools are
the Rellich theorm and the Lax-Milgram lemma.
The well-known theorem of Hilbert (see, e.g. \cite{CouHil}) and
regularity results like \cite[Prop.~I.2.2]{Temam}
lead to more precise results, namely:
\begin{corollary}
\label{STOeiFU}
The Stokes operator is an operator with a point spectrum only.
 All eigenvalues $\lambda_{j}$
of  ${\boldsymbol S}$ are real and of finite multiplicity.
The associated eigenfunctions 
$\{{\underline{w}}_{j}({\underline{x}})\}_{j=1}^{\infty}$
(counted in multiplicity)
are an orthogonal basis of
${\underline{\bb H}}(.)$ and ${\underline{\bb V}}(.)$:
\begin{align*}
{{(a)}}&\quad {\boldsymbol 
S}{\underline{w}}_{j}:=\lambda_{j}{\underline{w}}_{j}\quad\mbox{for
}\quad{\underline{w}}_{j}\in D({\boldsymbol S})
\quad\forall\,j=1,2,\dots\,;\\
{ (b)}& \quad
0\,<\lambda_{1}\leq\,\lambda_{2}\,\leq\cdots\leq\,\lambda_{j}
\,\leq\cdots\quad\mbox{and}\quad 
\lim_{j\rightarrow\infty}\lambda_{j}=\infty\,;
\\[-2mm]
{(c)}&\quad
\|{\underline{w}}_{j}\|_{{\underline{\bb 
H}}}\,=1\quad\forall \, j=1,2,\dots\,.
\end{align*}
\end{corollary}
\subsection{Eigenvalues and Eigenfunctions\label{sec_eigv}}
\noindent 
Formulas for  the
complete sets of Laplace and  
Stokes eigenfunctions on the circular annuli ${\Omega}^{*}_{\sigma}$ 
and on the unit ball ${\Omega}^{*}$ are derived in \cite{RumKug1} and
\cite{RumTh2024} as
(squares of) solutions of certain transcendental
equations. The
transformation to the ${\Omega}_{{\cal A}}$-domains then works as in
\cite{RuRuTh2016}. 
The ideas will be sketched below.
It ist worth to note, that -- like for the Bessel functions of the first kind and
orders  -- one can prove that
the consecutive zeros of our transcendental equations (taken for the
orders $k$ and $k+1$;  $k\,\in\, \{\frac{1}{2}+m,\,m\,\in\, {\bb N}_{o}\} $) 
separate interdependently. Also the first positive zero at order
$k$ is simple and smaller than 
the first positive zero at order $k+1$. The proof is similar to
the proof for the Bessel  functions 
(cf. \cite{AAR}), i.e. it uses standard tools such as the mean value theorem and
properties of Bessel functions. 
\pagebreak

\noindent For arbitrary ${\cal A}\,\in\,(0,\infty)$ the first simple eigenvalue
${\lambda}_{1,L}({\cal A})$ of the 
Laplacian (resp.~the first triple eigenvalue ${\lambda}_{1,S}({\cal A})$ of  the Stokes operator)
on ${\Omega}_{{\cal A}}$
are the squares of the smallest positive solutions
$\kappa_{1,L}({\cal A})$ and 
$\kappa_{1,S}({\cal A})$, respectively,
of the transcendental equations (see \cite{RumTh2024})
\begin{align}
0&=\,J_{\frac{1}{2}}(\kappa_{L}({\cal A}){(1+{\cal
   A}/{2}}))J_{-\frac{1}{2}}(\kappa_{L}({\cal A}){{\cal A}/{2}}) 
\,-\,J_{\frac{1}{2}}(\kappa_{L}({\cal A}){{\cal
   A}/{2}})J_{-\frac{1}{2}}(\kappa_{L}({\cal A}){(1+{\cal A}/{2}}))\,=\,\nonumber
\\
\label{trans1}
&=\,\frac{2}{\pi \kappa_{L}({\cal A})} [({\cal A}/{2})(1+{\cal A}/{2})]^{-\frac{1}{2}}
\sin(\kappa_{L}({\cal A}))   , 
\\
\label{trans2}
0&=\,J_{\frac{3}{2}}(\kappa_{S}({\cal A})({1+{\cal
   A}/{2
   }}))J_{-\frac{3}{2}}(\kappa_{S}({\cal A}){{\cal A}/{2}}) 
\,-\,J_{\frac{3}{2}}(\kappa_{S}({\cal A}){{\cal
   A}/{2}})J_{-\frac{3}{2}}(\kappa_{S}({\cal A})({1+{\cal A}/{2}}))\,. 
\end{align}
After the transformation to ${\sigma}$ we find the equivalent zeros
$\kappa_{1,L}({\sigma})$ and  
$\kappa_{1,S}({\sigma})$ as
the smallest positive solutions of 
\begin{align}
0&=\,J_{\frac{1}{2}}(\kappa_{L}({\sigma}))J_{-\frac{1}{2}}(\kappa_{L}({\sigma}){\sigma})
\,-\,J_{\frac{1}{2}}(\kappa_{L}({\sigma}){\sigma})J_{-\frac{1}{2}}(\kappa_{L}({\sigma}))\,=\,\nonumber
\\
\label{trans1sigma}
&=\,\frac{2}{\pi \kappa_{L}({\sigma})\sqrt{\sigma}}
\sin((1-{\sigma})\kappa_{L}({\sigma}))   , 
\\
\label{trans2sigma}
0&=\,J_{\frac{3}{2}}(\kappa_{S}({\sigma}))J_{-\frac{3}{2}}(\kappa_{S}({\sigma}){\sigma})
\,-\,J_{\frac{3}{2}}(\kappa_{S}({\sigma}){\sigma})J_{-\frac{3}{2}}(\kappa_{S}({\sigma}))\,.
\end{align}
The conversion formulas between the zeros are obvious for all ${\cal
  A}\,\in\,(0,\infty)$ and we see for the smallest ones that
\[
\kappa_{1,L}({\sigma})=({1+{\cal A}/{2}}) \kappa_{1,L}({\cal
  A})\,,
\qquad 
\kappa_{1,S}({\sigma})  =({1+{\cal A}/{2}})
\kappa_{1,S}({\cal A})\,.
\]
We note, that the formulas (\ref{trans1}) and  (\ref{trans1sigma})
are a simple consequence of the representation of spherical  Bessel  functions
by sinus and cosinus, see \eqref{Bess_Halbe} in the
Appendix. Now we choose the notation
$$
{\underline{\mathfrak{w}}}_{{\;\!}0}\,=\,\sin(\vartheta){\underline{\mathfrak{e}}}_{{\;\!}\varphi}\,, \quad{\underline{\mathfrak{w}}}_{{\;\!}-1}\,=\,cos(\varphi){\underline{\mathfrak{e}}}_{{\;\!}\vartheta}
-\sin(\varphi)\cos(\vartheta){\underline{\mathfrak{e}}}_{{\;\!}\varphi}\,,\quad {\mbox{and}} \,
\quad
{\underline{\mathfrak{w}}}_{{\;\!}1}\,=\,sin(\varphi){\underline{\mathfrak{e}}}_{{\;\!}\vartheta}
+\cos(\varphi)\cos(\vartheta){\underline{\mathfrak{e}}}_{{\;\!}\varphi}\,.\nonumber
$$
The eigenfunctions ${{w}}_{{\;\!}1,L,{\cal A}}({\underline{x}})$ and 
${\underline{w}}_{{\;\!}1,\alpha,S,{\cal A}}({\underline{x}})$ in
${\Omega}_{{\cal A}}$ (here $\alpha\in\{-1,0,1\})$ are 
\begin{align}
{{w}}_{{\;\!}1,L,{\cal A}}({\underline{x}})&=\frac{\tilde{c}_{1,L,{\cal A}}}{\sqrt{r}}\left(
J_{\frac{1}{2}}(\pi r)\,-\,
\frac{J_{{\frac{1}{2}}}({\textstyle{\frac{ \pi {\cal A}}{2}}})}{J_{-{\frac{1}{2}}}({\textstyle{\frac{ \pi {\cal A}}{2}}})}
J_{-\frac{1}{2}}(\pi r )\right) \label{eiflaplace} \quad\quad \text{and}
\\
{\underline{w}}_{{\;\!}1,\alpha,S,{\cal A}}({\underline{x}})&=\frac{\tilde{c}_{1,\alpha,S,{\cal A}}}{\sqrt{r}}
\left(
J_{\frac{3}{2}}(\kappa_{1,S}({\cal A})r)\,-\,
\frac{J_{{\frac{3}{2}}}({\textstyle{\frac{ \kappa_{1,S}({\cal A}) {\cal A}}{2}}})}
{J_{-{\frac{3}{2}}}({\textstyle{\frac{ \kappa_{1,S}({\cal A}) {\cal A}}{2}}})}
J_{-\frac{3}{2}}( \kappa_{1,S}({\cal A}) r
)
\right)
                                                              {\underline{\mathfrak{w}}}_{{\;\!}\alpha}(\vartheta,\varphi)\,.
                                                              \label{eifstokes}
\end{align}
The numbers $\tilde{c}_{1,L,{\cal A}}$ and $\tilde{c}_{1,\alpha,S,{\cal A}}$ are
scaling (to $1$) constants in the ${{\bb L}}_{{2}}(.)$- 
resp. ${\underline{\bb L}}_{{\;\!}{2}}(.)$-sense (see also
\cite{RumTh2024}). The notations of the 
${\underline{\mathfrak{w}}}_{{\;\!}\alpha}$ 
refer to the spherical surface harmonics
$Z_{1}^{\alpha}(.) $. Here $\alpha$ counts as above, i.e. $\alpha\in\{-1,0,1\}$. 
\begin{corollary}
\label{PiQuadrat}
For all ${\cal A}\,\in\,(0,\infty)$ the first simple eigenvalue ${\lambda}_{1,L}({\cal A})$ of the 
Laplacian ${\boldsymbol L}$ 
is ${\lambda}_{1,L}({\cal A})\,=\,{\pi}^2$.
\end{corollary}
\begin{proof}
We use  Eq.  \eqref{trans1}, more precisely, $
0\,=\,\frac{2}{\pi \kappa_{L}({\cal A})} (({\cal A}/{2})(1+{\cal A}/{2}))^{-\frac{1}{2}}
\sin(\kappa_{L}({\cal A})) \,\,.$
For all 
${\cal A}\,\in\,(0,\infty)$ its smallest positive solution is
$\kappa_{1,L}({\cal A})\,=\,{\pi}$.  
Finally we conclude that ${\lambda}_{1,L}({\cal A})\,=\,\kappa_{1,L}^2({\cal A})
\,=\,{\pi}^2$.
\end{proof}
\noindent {\bf Remark:} 
We know for the one-dimensional Laplace problem that the eigenvalue ${\lambda}_{1}\,=\,{\pi}^2\,$.
This is also the idea behind the small gap limit in three
dimensions. Namely, 
for the first simple eigenvalue it holds for ${\cal
  A}\,\to\,\infty$ that ${\lambda}_{1,L}({\cal A})\,\to\,{\pi}^2$. 
In particular, we can take 
\eqref{eq:poincarepolarradial} (cf. also \eqref{eq:poincareEins} and \eqref{eq:poincareZwei})
for $\tilde w$ varying in 
$ {\bb W}_{2}^{1}\hspace{-.62cm}{~}^{{~}^{{~}^{o}}}\hspace{.2cm}(\frac{\mathcal A}{2},1+\frac{\mathcal A}{2})$ to get
for all ${\cal A}\,\in\,(0,\infty)$ with \cite{nazarov2000}  that 
\begin{equation}
\label{hallos}
\left(
c_{p}({\cal A})
\right)^2 
=
\left(\frac{1}{\pi}\right)^2 = 
\max_{\tilde w } \frac{
\int_{\frac{\mathcal A}{2}}^{1+\frac{\mathcal A}{2}}
\abs{\tilde w(r)}^2 \, dr}
{
\int_{\frac{\mathcal A}{2}}^{1+\frac{\mathcal A}{2}}
\abs{\tilde w'(r)}^2 \, dr}\,.
\end{equation}
\section{Investigation of the Limiting Cases \label{Inv_Lim}}
Subsequently we will study the limiting cases ${\cal A}\,\to\,0$
and ${\cal A}\,\to\,{\infty}$  
separately, because the methods are completely different.
It is quite obvious that the cases ${\cal A}\,\to\,0$ and  ${\sigma}\,\to\,0$ are
identical and we formally obtain the punctured ball
\[
{\Omega}^{*}_{\setminus \{{\underline{0}}\}} :=  
{\Omega}^{*}\setminus \{{\underline{0}}\} \quad\mbox{in the limits}\quad
 \lim_
{{\sigma}\,\to\,0}{\Omega}^{*}_{\sigma}={\Omega}^{*}_{\setminus
  \{{\underline{0}}\}}  \,,\quad 
 \lim_{{\cal A}\,\to\,0}{\Omega}_{\cal A}={\Omega}^{*}_{\setminus
  \{{\underline{0}}\}} \,,\quad 
{\Omega}^{*}_{\setminus
  \{{\underline{0}}\}}  \simeq {\Omega}^{*}\,.
\]
For simplicity we will choose the parameter ${\sigma}$ and may use
(\ref{calA}) to connect ${\cal A}$ and 
${\sigma}$ for the (first) eigenvalues. So we note, e.g., for the
smallest eigenvalue (the square of the smallest positive zero of
\eqref{trans1}) that 
\begin{align}
\lambda_{1,L}({\sigma})\,=\,
(\kappa_{1,L}({\sigma}))^2\,=\,(({1+{\cal A}/{2}}) \kappa_{1,L}({\cal A}))^2\,=\,({1+{\cal A}/{2}})^2 \lambda_{1,L}({\cal A})\,.
\end{align}
The parameter 
${\cal A}$ is replaced by the radius $R\,={\cal A}/2$ for
convenience. Looking for the 
limit of 
${\Omega}_{{\cal A}}$ as ${\cal A}\to \infty$ one easily sees, that
they are ``losing their curvature'' and in the limit become an
unbounded layer also known as the {\em small gap limit}
\[
{\Omega}_{{\cal A} =\infty}\,:=
\{{(s,\tau,t)^{T}} \in {\bb R}^{3}:\,0 < s < 1\}\,.
\]
It is worth to note,
that in the limiting process ${\cal A}\,\to\,{\infty}$ (or $R\,={\cal
  A}/2\to\infty$) the Laplace and 
Stokes operator lose the property of being operators with a pure point spectrum. 
But nevertheless, they stay linear positive operators and this ensures
the existence of smallest eigenvalues. 
In contrast to our ideas for coordinate transformations in
\cite{RuRuTh2016} we are going to employ 
considerations of the limiting value. There we will use the
transcendental equations for the eigenvalues 
in the limiting process
${\cal A}\,\to\,{\infty}$ (or $R\,={\cal A}/2\,\to\,{\infty}$). 
Let us repeat the transcendental equations
\begin{align}
0&=\,J_{\frac{1}{2}}(\kappa_{L}({\cal A}){(1+{\cal
   A}/{2}}))J_{-\frac{1}{2}}(\kappa_{L}({\cal A}){{\cal A}/{2}}) 
\,-\,J_{\frac{1}{2}}(\kappa_{L}({\cal A}){{\cal
   A}/{2}})J_{-\frac{1}{2}}(\kappa_{L}({\cal A}){(1+{\cal A}/{2}}))
   \label{transoo}
  \\
\mbox{and}\qquad 0&=\,J_{\frac{3}{2}}(\kappa_{S}({\cal A})({1+{\cal
   A}/{2
   }}))J_{-\frac{3}{2}}(\kappa_{S}({\cal A}){{\cal A}/{2}}) 
\,-\,J_{\frac{3}{2}}(\kappa_{S}({\cal A}){{\cal
   A}/{2}})J_{-\frac{3}{2}}(\kappa_{S}({\cal A})({1+{\cal A}/{2}}))\,\label{transoo2}. 
\end{align}
\label{sec_theo}
\subsection{The behaviour of the Laplace eigenvalues for
  $\boldsymbol{\sigma\,\to\,0}$  
\label{Green} }
We use the notations ${\boldsymbol L}_{\sigma}$ for the  Laplace operator on
${\Omega}^{*}_{\sigma}$ and ${\boldsymbol L}_{\sigma}^{-1}$ for its inverse at ${\sigma}\,\in\,(0,1)$.
Additionally we denote the Laplace 
operator on the unit ball
${\Omega}^{*}$ by ${\boldsymbol L}_{o}$ and its inverse 
${\boldsymbol L}_{o}^{-1}$, respectively.
Let us use the Green's functions for the
(negative) Laplacians to define the inverse operators: 
Let ${\underline{x}}$ and ${\underline{y}}$ be two
points in 
${\overline{{\Omega}}}^{{\;\!}*}$. 
Then Green's function for the
(negative) Laplacian on ${\Omega}^{*}$ with zero trace is given by (cf. \cite{CouHil}):
\begin{equation*}
G({\underline{x}},{\underline{y}})=\frac{1}{4\pi} \left(\frac{\displaystyle{1}}{\displaystyle{
\|{\underline{x}}-{\underline{y}}\|}}\,-\,
\frac{\displaystyle{1}} 
{\displaystyle{(\|{\underline{x}}\|^{2}
\|{\underline{y}}\|^{2}+1-2\|{\underline{x}}\|
\|{\underline{y}}\|\cos({\underline{x}},{\underline{y}}))^{\frac{1}{2}}}}\,\right).
\end{equation*} 
The symmetric form of the Green's function for the (negative)
Laplacian with vanishing traces on $\partial{\Omega}^{*}_{\sigma}=
\omega_{\sigma} \cup \omega$ on 
every ${\Omega}_{\sigma}^{*}$ (for fixed $\sigma\,>\,0$) is (cf. \cite{GilTru})
\begin{eqnarray}\label{sternstern}
G_{\sigma}({\underline{x}},{\underline{y}})=G({\underline{x}},{\underline{y}})+
\sum_{k\in {\bb N}}\frac{{\sigma}^{k}}{4\pi}\bigg(\frac{1}{\displaystyle{\|{\underline{y}}-{\sigma}^{2k}{\underline{x}}\|}}
+\frac{1}{\displaystyle{\|{\sigma}^{2k}{\underline{y}}-{\underline{x}}\|}}+
\hspace*{2cm}\quad\quad\nonumber\\ 
\frac{-1}
{\displaystyle{(\|{\underline{x}}\|^{2}
\|{\underline{y}}\|^{2}+{\sigma}^{4k}-2{\sigma}^{2k}\|{\underline{x}}\|
\|{\underline{y}}\| 
\cos({\underline{x}},{\underline{y}}))^{\frac{1}{2}}}}+\frac{-1}
{\displaystyle{({\sigma}^{4k}\|{\underline{x}}\|^{2}
\|{\underline{y}}\|^{2}+1-2{\sigma}^{2k}\|{\underline{x}}\|
\|{\underline{y}}\| 
\cos({\underline{x}},{\underline{y}}))^{\frac{1}{2}}}}\bigg)\,\,,
\end{eqnarray}
where the points ${\underline{x}}$ and ${\underline{y}}$ have to be in
${\overline{{\Omega}}}^{{\;\!}*}_{\sigma}$. This shape of the Green's function is useful for the study 
of ${\sigma}\,\to\,0$. By a straightforward calculation the limit becomes
\begin{align}
\label{Lapl_GWsig}
\lim_{{\sigma}\,\to\,0}{\;\!}
  G_{\sigma}({\underline{x}},{\underline{y}})=G({\underline{x}},{\underline{y}})
  \,. 
\end{align}
The basic properties of the eigenfunctions of  ${\boldsymbol L}_{o}$  
and Bessel's differential operator can be found in \cite[Chapters 5
and 8]{Triebel} together with the explicit 
representation of the inverse Laplacian by  Fourier series in 
the eigenfunctions. \\[.1cm]
The first eigenvalue   $\lambda_{1,L}({o})=
(\kappa_{1,L}({o}))^{2}$ 
of  the operator ${\boldsymbol L}_{o}$ is simple, where $\kappa_{1,L}({o})\,=\,\pi$ is the  
first positive root of $J_{\frac{1}{2}}(.)$.
The corresponding first 
eigenfunction is
$w_{1,o}({\underline{x}})=\|{\underline{x}}\|^{-1}c_{1,o}J_{\frac{1}{2}}(\pi\|{\underline{x}}\|)$
with constant $c_{1,o}$ such that $\| w_{1,o}({\underline{x}})
\|_{{\bb L}_{2}(.)}=1$. \\
The operators 
${\boldsymbol L}_{\sigma}$  have the same property: Their first eigenvalues 
$\lambda_{1,L}({\sigma})\,=\,(1-\sigma)^{-2}{\pi}^{2}$ are simple and the respective 
$w_{1,{\sigma}}({\underline{x}})$ are their first eigenfunctions.
\\[3mm]
The operators ${\boldsymbol L}_{o}^{-1}$ and ${\boldsymbol
  L}_{\sigma}^{-1}$ are selfadjoint, positive and 
compact. Thus, we may write:
\begin{align}
\label{neu}
    \pi^{-2}\,=\, (\lambda_{1,L}({o}))^{-1}&= \max_{u \in {\bb
        L}_{2}({\Omega}^{*}_{{~}}): \| u \|_{{\bb L}_{2}(.)}=1 }{\;\!}
    \int_ {{\Omega}^{*}_{{~}}} \int_ {{\Omega}^{*}_{{~}}}
    G({\underline{x}},{\underline{y}})u({\underline{y}})
    u({\underline{x}})d{\underline{y}}d{\underline{x}}=\|{\bf
      L}_{o}^{-1}\|
    \\
    (1-\sigma)^{2}{\pi}^{-2} \,=\, (\lambda_{1,L}({\sigma}))^{-1}&=
    \max_{u \in {\bb L}_{2}({\Omega}^{*}_{\sigma}): \| u \|_{{\bb
          L}_{2}(.)}=1 }{\;\!}  \int_ {{\Omega}^{*}_{\sigma}} \int_
    {{\Omega}^{*}_{\sigma}}
    G_{\sigma}({\underline{x}},{\underline{y}})
    u({\underline{y}})u({\underline{x}})d{\underline{y}}d{\underline{x}}=\|{\bf
      L}_{\sigma}^{-1}\|\,\nonumber,
\end{align}
where the maxima are attained at $u=w_{1,o}({\underline{x}})$ and at
$u=w_{1,{\sigma}}({\underline{x}})$, 
respectively. The equations above immediately show the behaviour of the Laplace eigenvalues for
${\sigma\,\to\,0}$. We note, that the theoretical background is a
simple tool of classical partial differential equations.
We sketch the ideas by using again
$
(\lambda_{1,L}({o}))^{-1}=\|{\bf L}_{o}^{-1}\|
\quad \mbox{and}\quad
(\lambda_{1,L}({\sigma}))^{-1}=\|{\bf L}_{\sigma}^{-1}\|\,. $
For ${\sigma}\in\,(0,1)$ and $ {\sigma}_{o}\geq {\sigma}$ obviously, 
${\Omega}_{{\sigma}_{o}}^{*}\subset {\Omega}_{{\sigma}}^{*}\subset 
{\Omega}^{*}$. Thus, the values
$(\lambda_{1,L}({\sigma}))^{-1}$ constitute a continuous, 
monotonically decreasing function 
with respect to ${\sigma}$ (applying
\cite[Satz 3 on page 355 in Bd. 1]{CouHil}).
Also for any fixed ${\sigma}_{o}<1$ we see 
\begin{equation}\label{stern}
(\lambda_{1,L}({\sigma}_{o}))^{-1}\leq  (\lambda_{1,L}({\sigma}))^{-1}
\leq (\lambda_{1,L}({o}))^{-1}
\quad
\mbox{for}\quad {\sigma}\in\,(0,{\sigma}_{o})\,. 
\end{equation}
The second tool for our proof is the application of {\em constrained
  subsets} of ${\bb L}_{2}({\Omega}^{*})$. 
Let 
\[
{\bb L}_{2}({\Omega}^{*}_{\backslash {\sigma}}):=\{v\,\in\,{\bb
  L}_{2}({\Omega}^{*}): v=0\,\, {\mbox{a.e. on }} {\Omega}^{*} 
\setminus {\Omega}^{*}_{\sigma}\}
\]
be the space of almost everywhere vanishing functions for 
${\underline{x}} \in {\bb R}^{3}$ with $\|{\underline{x}}\|\,<{\sigma}$. It
is obvious, that 
${\bb L}_{2}({\Omega}^{*}_{\backslash {\sigma}})\simeq {\bb
  L}_{2}({\Omega}^{*}_{\sigma})$ 
in the sense of an identity map.
We apply the standard property of maxima on subsets combined with the
central tool of the limiting process  
(\ref{Lapl_GWsig}) to see, that there exists 
\[ 
\lim_{{\sigma}\,\to\,0}{\;\!}  (\lambda_{1,L}({\sigma}))^{-1}
\leq  (\lambda_{1,L}({o}))^{-1}
\]
and $\lambda_{1,L}({\sigma}))$ is continuous in $0$ (see \cite{Weid}).
Using that for every $u$
${\bb L}_{2}({\Omega}^{*}_{\setminus \{{\underline{0}}\}} )$ we also
have $u \in {\bb L}_{2}({\Omega}^{*} )$ and the definitions of $G $
and $G_\sigma$ we conclude
\begin{align}\label{neu1}
\int_{ {\Omega}^{*}_{{~}}}\int_ {{\Omega}^{*}_{{~}}} G({\underline{x}},{\underline{y}})
u({\underline{y}})u({\underline{x}})d{\underline{y}}d{\underline{x}}=
\int_ {{\Omega}^{*}_{\setminus \{{\underline{0}}\}}}\int_ {{\Omega}^{*}_{\setminus \{{\underline{0}}\}}}
\lim_{{\sigma}\,\to\,0} G_{\sigma}({\underline{x}},{\underline{y}}) 
u({\underline{y}})u({\underline{x}})d{\underline{y}}d{\underline{x}}\quad\forall
  u \in {\bb L}_{2}({\Omega}^{*} )\,\,\,. 
\end{align}
From \eqref{neu}, \eqref{neu1} and the continuity of $\lambda_{1,L}({\sigma}))$ in $0$ we may conclude that 
\begin{align*}
\max_{u \in {\bb L}_{2}({\Omega}^{*}_{{~}}):  \| u \|_{{\bb L}_{2}(.)}=1  }{\;\!}&
\int_ {{\Omega}^{*}_{\setminus \{{\underline{0}}\}}}\int_ {{\Omega}^{*}_{\setminus \{{\underline{0}}\}}}
\lim_{{\sigma}\,\to\,0} G_{\sigma}({\underline{x}},{\underline{y}}))
u({\underline{y}})u({\underline{x}})d{\underline{y}}d{\underline{x}}
                           \\ 
=&\max_{u \in {\bb L}_{2}({\Omega}^{*}_{{~}}):   \| u \|_{{\bb L}_{2}(.)}=1  }{\;\!}
\int_{ {\Omega}^{*}_{{~}}}\int_ {{\Omega}^{*}_{{~}}} G({\underline{x}},{\underline{y}})
u({\underline{y}})u({\underline{x}})d{\underline{y}}d{\underline{x}}\quad
\\
  =&(\lambda_{1,L}({0}))^{-1}\\
  =&\,\,\lim_{{\sigma}\,\to\,0} (\lambda_{1,L}({\sigma}))^{-1}=
\lim_{{\sigma}\,\to\,0}
\Big( \max_{u \in {\bb L}_{2}({\Omega}^{*}_{\sigma}):   \| u \|_{{\bb L}_{2}(.)}=1  }{\;\!}
\int_ {{\Omega}^{*}_{\sigma}}   \int_ {{\Omega}^{*}_{\sigma}}  (
G_{\sigma}({\underline{x}},{\underline{y}}))
u({\underline{y}})u({\underline{x}})d{\underline{y}}d{\underline{x}}
\Big).
\end{align*}
We note, that also the limit in the eigenfunctions 
\begin{align*}
w_{1,o}({\underline{x}})=\lim_{{\sigma}\,\to\,0}{\;\!}w_{1,{\sigma}}({\underline{x}})
\end{align*}
is well-defined in the sense of almost uniform convergence. Here
$w_{1,{\sigma}}({\underline{x}})$ is regarded as the continuous 
extension of $w_{1,{\sigma}}({\underline{x}})$ on
${\overline{{\Omega}}}^{{\;\!}*}$ by
zero in $\,{\overline{{\Omega}}}^{{\;\!}*}
\setminus{\overline{{\Omega}}}^{{\;\!}*}_{\sigma} $.
It is worth noticing, that the behaviour of the Green's functions for the
(negative) Laplacians in the process ${\sigma\,\to\,0}$ may be seen as 
forgetting the point ${\underline{0}}$ (with its zero boundary
conditions) in the limit. 
\subsection{The behaviour of the first Stokes eigenvalues for
  $\boldsymbol{\sigma\,\to\,0}$ \label{StoJ1Y1} } 
We denote the  Stokes operators on
 ${\Omega}^{*}_{\sigma}$ by ${\bf S}_{\sigma}$  and the  ones 
on the unit ball by ${\bf S}_{o}$ (as in Subsection~\ref{Green}).
The elements of the family 
${\bf S}_{\sigma}$
 are selfadjoint and positive just as ${\bf S}_{o}$.
For the first triple eigenvalues $\lambda_{1,S}({o})$ and
 $\lambda_{1,S}({\sigma})$ we observe 
 \begin{align}
\label{Stok_Eig}
\lambda_{1,S}({o})=
\min_{{\underline{u}} \in {\underline{\bb S}}_{{\;\!}}^{2}(({\Omega}^{*}_{{~}}):
   \| {\underline{u}} \|_{{\underline{\bb H}}_{{\;\!}}(.)
   }} ({\bf
   S}_{o}{\underline{u}},{\underline{u}})_{{\underline{\bb
   H}}_{{\;\!}}(.)
   } 
\quad 
{\mbox{and}}\quad
\lambda_{1,S}({\sigma})=
\min_{{\underline{u}} \in {\underline{\bb S}}_{{\;\!}}^{2}(({\Omega}^{*}_{{\sigma}}):
\| {\underline{u}} \|_{{\underline{\bb H}}_{{\;\!}}(.) }} 
({\bf   S}_{{\sigma}} {\underline{u}},{\underline{u}})_{{\underline{\bb H}}_{{\;\!}}(.)  }
\,, 
\end{align}
where the minima are attained at
${\underline{u}}\in{\mbox{span}}\{{\underline{w}}_{{\;\!}1,\alpha,o}({\underline{x}})\}_{\alpha = -1}^{1}$ and at  
${\underline{u}}\in{\mbox{span}}\{{\underline{w}}_{{\;\!}1,\alpha,\sigma}({\underline{x}})\}_{\alpha = -1}^{1}$,
respectively. 
For ${\sigma}\in\,(0,1)$ the eigenvalues
$\lambda_{1,S}({\sigma})$ constitute a continuous, monotonically
increasing function of ${\sigma}$. This is due to \cite[Satz 3 on p.~355 in Bd.~1]{CouHil}
and the inclusions  
${\Omega}_{{\sigma}_{o}}^{*}\subset {\Omega}_{{\sigma}}^{*}\subset 
{\Omega}^{*}$ for any $ {\sigma}_{o}\geq {\sigma}$. 
In addition we
observe 
$\lambda_{1,S}({\sigma}_{o})\geq  \lambda_{1,S}({\sigma}) \geq
\lambda_{1,S}({o})$ 
for ${\sigma}\in\,(0,{\sigma}_{o})$, for any fixed ${\sigma}_{o}<1$. 
With similar arguments as in Subsection \ref{Green} we
see that the right-hand limit  
of the first Stokes eigenvalues attains the unique value
 \begin{align}\label{eWSlimsigma_0}
  (\kappa_{1,S}({o}))^{2}\,\equiv \,
 {\lambda_{1,S}({o})}=\lim_{{\sigma}\,\to\,0}{\;\!}
   {\lambda_{1,S}({\sigma})}\,\equiv \, 
\lim_{{\sigma}\,\to\,0}{\;\!}  \big( {\kappa_{1,S}({\sigma})}\big)^{2}\,.
 \end{align}
 One achieves this result by Eq. 
 (\ref{Stok_Eig}).  We use the formal construction with {\em constrained
  subsets} of ${\underline{\bb H}}_{{\;\!}}({\Omega}^{*})$ to
define the space of almost everywhere vanishing vector functions in 
${\underline{x}} \in {\bb R}^{3}:\|{\underline{x}}\|\,<{\sigma}$. 
\[{\underline{\bb H}}_{{\;\!}}({\Omega}^{*}_{\backslash {\sigma}}):=\{\underline{v}\,\in\,{\underline{\bb H}}_{{\;\!}}({\Omega}^{*}): \underline{v}=\underline{0}\,\, 
{\mbox{a.e. on }} {\Omega}^{*} 
\setminus {\Omega}^{*}_{\sigma}\}\,.
\]
Like in Subsection~\ref{Green} we find that 
${\underline{\bb H}}_{{\;\!}}({\Omega}^{*}_{\backslash {\sigma}})\simeq {\underline{\bb H}}_{{\;\!}}({\Omega}^{*}_{\sigma})$ 
in the sense of an identity map. Additionally,  ${\underline{\bb
    H}}_{{\;\!}}({\Omega}^{*}_{\setminus \{{\underline{0}}\}})$ is 
equivalent to 
${\underline{\bb H}}_{{\;\!}}({\Omega}^{*})$ in the sense of an identity map, where
${\Omega}^{*}_{\setminus \{{\underline{0}}\}}$ was explained in Section~\ref{Inv_Lim}.
Now we use the monotonicity in ${\sigma}$ and 
polar coordinates (with $r=\|{\underline{x}}\|$)
(cf. \eqref{eifstokes}) and see that
\begin{align}
\label{efSsigma_0}
{\underline{w}}_{{\;\!}1,\alpha,o}({\underline{x}})&=
\frac{\tilde{c}_{1,\alpha,o}}{\sqrt{r}}
J_{{\frac{3}{2}}}(\kappa_{1,S}({o})r)\,{\underline{\mathfrak{w}}}_{{\;\!}\alpha}(\vartheta,\varphi)
\,
&& {\mbox{in}} \,\, {\Omega}^{*} \nonumber\\
{\underline{w}}_{{\;\!}1,\alpha,\sigma}({\underline{x}})&=\frac{\tilde{c}_{1,\alpha,\sigma}}{\sqrt{r}}
\Big(
J_{\frac{3}{2}}(\kappa_{1,S}({\sigma})r)\,-\,
\frac{J_{{\frac{3}{2}}}({\textstyle{ \kappa_{1,S}({\sigma}) {\sigma}}})}
{J_{-{\frac{3}{2}}}({\textstyle{{ \kappa_{1,S}({\sigma}) {\sigma}}}})}
J_{-\frac{3}{2}}( \kappa_{1,S}({\sigma}) r
)
\Big)
{\underline{\mathfrak{w}}}_{{\;\!}\alpha}(\vartheta,\varphi)
\quad
&& {\mbox{in}} \,\, {\Omega}_{{\sigma}}^{*}.
\end{align} 
Here $\tilde{c}_{1,\alpha,o}$ and $\tilde{c}_{1,\alpha,\sigma}$  $(\alpha\in\{-1,\,0,\,1\}$
are scaling (to $1$)
constants in the ${\underline{\bb L}}_{{\;\!}{2}}(.)$-sense. 
We note, that the boundedness of these constants is obvious, but their
regularity as functions  
of ${\sigma}$ will be clarified only later on.
We omit these constants at first.
For the behaviour of the zeros ${\kappa_{1,S}({\sigma})}$ of the
transcendental equations and of the zeros of  
Besselfunctions ${\kappa_{1,S}({o})}$ we refer to
Subsection~{\ref{sec_eigv}}.
\\[3mm]
We investigate the limit ${{\sigma}\,\to\,0}$ term by term
regarding the functions
${\underline{w}}_{{\;\!}1,\alpha,{\sigma}}({\underline{x}})$  as the
continuous  
extension of ${\underline{w}}_{{\;\!}1,\alpha,{\sigma}}({\underline{x}})$  onto
${\overline{{\Omega}}}^{{\;\!}*}$ by  $\underline{0}$
(the zero vector) in 
${\overline{{\Omega}}}^{{\;\!}*}
\setminus{\overline{{\Omega}}}^{{\;\!}*}_{\sigma} $ 
if necessary. 
It follows from (\ref{eWSlimsigma_0}), that for all
${\underline{x}}\,\in\,{\overline{{\Omega}}}^{{\;\!}*}$
($r\,\in\,[0,1]$) 
\begin{align}\label{efGW1Ssigma_0}
\lim_{{\sigma}\,\to\,0}{\;\!} \frac{J_{\frac{3}{2}}(\kappa_{1,S}({\sigma})r)}{\sqrt{r}}=
\frac{J_{\frac{3}{2}}(\kappa_{1,S}({o})r)}{\sqrt{r}}\,.
\end{align} 
Furthermore, using the properties of the Bessel function
${J_{-\frac{3}{2}}(.)}$  
with \cite{nazarov2000} we find
for all $r,{\sigma}\,\in\,[0,1]$ with $ r \geq
{\sigma}$ that
\begin{align}\label{Bessel1frac}
\left\|
  \frac{J_{-\frac{3}{2}}(\kappa_{1,S}({\sigma})r)}{J_{-\frac{3}{2}}
(\kappa_{1,S}({\sigma}){\sigma})}\right\|
  \,\leq\,1= 
\lim_{r,{\sigma}\,\to\,0: \,r \geq {\sigma}}{\;\!} \left\|\frac{J_{-\frac{3}{2}}
(\kappa_{1,S}({\sigma})r)}{J_{-\frac{3}{2}}(\kappa_{1,S}({\sigma}){\sigma})}\right\|\,.
\end{align} 
Finally, from 
\begin{align*}
0\,\leq\,
\lim_{{\sigma}\,\to\,0}{\;\!} \frac{J_{\frac{3}{2}}(\kappa_{1,S}({\sigma})\sigma)}{\sqrt{r}}\,\leq\,
\lim_{{\sigma}\,\to\,0}{\;\!}
\frac{J_{\frac{3}{2}}(\kappa_{1,S}({{\sigma}}){\sigma})}{\sqrt{{\sigma}}}\,=\,0
\end{align*} 
and Eqs.~(\ref{efGW1Ssigma_0}) and  (\ref{Bessel1frac}) we may conclude that 
$\lim_{{\sigma}\,\to\,0}{\;\!}
{\underline{w}}_{{\;\!}1,\alpha,{\sigma}}({\underline{x}})=
{\underline{w}}_{{\;\!}1,\alpha,o}({\underline{x}})$ 
in the sense of pointwise convergence on
${\overline{{\Omega}}}^{{\;\!}*}$ and in the sense 
of almost uniform convergence too due to the obvious right-hand
continuity of $\tilde{c}_{1,\alpha,\sigma}$. 
\subsection{The behaviour of the first Laplace and Stokes 
eigenvalues for 
$\boldsymbol{{\cal A}\,\to\,{\infty}}$ \label{LaStoInf} } 
We directly investigate the transcendental equations for the first eigenvalues 
of the Laplace and the
Stokes operators as ${\cal A}\,\to\,{\infty}$.  	
In contrast to our 2d studies in \cite{RuRuTh2016} 
the key to success is not the use of a  
coordinate transformation like:
\begin{equation}
\label{Coordinate}
\mbox{Let }R:=\frac{\cal A}2\qquad\mbox{and}\quad
\left\{\begin{aligned}
r&:=R+s &&\mbox{for $s\,\in\,[0,1]$}\\
{\vartheta}&:={\tau}/{R }
&&\mbox{for } \tau\,\in\,[0,R\pi]\,\\
{\varphi}&:={t}/{R }
&&\mbox{for } t\,\in\,[-R\pi,R\pi]\,.
\end{aligned}\right.
\end{equation}
Only the transformation $r:=R+s$ for $s\,\in\,[0,1]$ is useful here. The so-called
small-gap limit is more or less
to establish a relationship between the first eigenvalues 
of the Laplace operator (resp. the
Stokes operator) and the one-dimensional Laplace eigenvalue problem as ${\cal A}\,\to\,{\infty}$ 
\cite[Prop.~1.1]{nazarov2000}, where 
$\tilde w$ varies  in $ {\bb
  W}_{2}^{1}\hspace{-.62cm}{~}^{{~}^{{~}^{o}}}\hspace{.2cm}(\frac{\mathcal
  A}{2},1+\frac{\mathcal A}{2})$ and 
\begin{equation}
\label{eq:1DEigenwert}
\hspace*{4cm}
{\pi}^2 = 
\min_{\tilde w } \frac{
\int_{\frac{\mathcal A}{2}}^{1+\frac{\mathcal A}{2}}
\abs{\tilde w'(r)}^2 \, dr}
{
\int_{\frac{\mathcal A}{2}}^{1+\frac{\mathcal A}{2}}
\abs{\tilde w(r)}^2 \, dr}\,\,.
\end{equation}
We are going to recover the number ${\pi}^2$ in the following Corollary:
\begin{corollary}
\label{PiQuadratinf}
For the first simple eigenvalue ${\lambda}_{1,L}({\cal A})$ of the 
Laplacian ${\boldsymbol L}$ 
we obtain that 
$\lim_{{\cal A}\to\infty} {\lambda}_{1,L}({\cal A})\,=\,{\pi}^2$.
\end{corollary}
\begin{proof} 
Again we use Eqs.~\eqref{trans1} and
\eqref{Bess_Halbe}. 
We know from Corollary \ref{PiQuadrat}, that for all 
${\cal A}\,\in\,(0,\infty)$ the smallest positive
solution of Eq.~\eqref{trans1} is  
$\kappa_{1,L}({\cal A})\,=\,{\pi}$ and that
$\kappa_{1,L}({0})\,=\,{\pi}$ (cf. Subsection \ref{Green}). 
In what follows we will use the abbreviation $\kappa
\,=\,\kappa_{L}({\cal A})$. From \eqref{trans1} and
\eqref{Bess_Halbe} we infer that
\begin{align}\label{transinf1}
\hspace*{-1.4cm}
0\,=\,J_{\frac{1}{2}}(\kappa{(1+{\cal
   A}/{2}}))J_{-\frac{1}{2}}(\kappa{{\cal A}/{2}}) 
\,-\,J_{\frac{1}{2}}(\kappa{{\cal
   A}/{2}})J_{-\frac{1}{2}}(\kappa{(1+{\cal A}/{2}}))\,=
   \,\frac{2}{\pi \kappa} [({\cal A}/{2})(1+{\cal A}/{2})]^{-\frac{1}{2}}
\sin(\kappa)\,\,.
\end{align}
Taking into consideration the continuity of \eqref{transinf1} in
$\kappa\,(=\,{\pi})$ this is equivalent to
\begin{align}
 0\,=\,  \,\frac{2}{\pi \kappa} 
\sin(\kappa) \,\,,
\end{align} 
with the smallest positive solution $\kappa\,=\,{\pi}$.
Finally we use the relation $\lim_{{\cal A}\to\infty} {\lambda}_{1,L}({\cal A})\,=\,
\lim_{{\cal A}\to\infty} \kappa^2_{1,L}({\cal A})
\,=\,{\pi}^2$.
\end{proof}
\noindent We can read the asymptotic behaviour
of $c_{p}({\cal A})$ as (cf. (4)).
\begin{align*}
\lim_{{\cal A}\,\to\,\infty}{\;\!} c_{p}({\cal A})={\frac{1}{\pi}}\,\,.
\end{align*}
The 'approximate' limit of the first eigenfunction $w_{1,L,{\cal
    A}}({\underline{x}})$ is now immediately obtained (cf.
    \eqref{Coordinate}) as
\begin{align*}
w_{1,L,{\cal A}}({\underline{x}})\,\approx\,{\tilde{c}}\,{\sin}(\pi s)
  \qquad{\mbox{for large}}\,{\cal A}\,.
\end{align*}
This limit can be calculated with \eqref{eiflaplace} and $r\,:=\,\frac{\cal A}{2}+s\,=\,R+s \,$
\,{for $s\,\in\,[0,1]$} :
\begin{align*}
w_{1,L,{\cal A}}({\underline{x}})\,=\,\frac{{\tilde{c}}_{1,L,{\cal A}}}{\sqrt{r}}\left(
J_{\frac{1}{2}}(\pi r)\,-\,
\frac{J_{{\frac{1}{2}}}({\textstyle{\frac{ \pi {\cal A}}{2}}})}{J_{-{\frac{1}{2}}}({\textstyle{\frac{ \pi {\cal A}}{2}}})}
J_{-\frac{1}{2}}(\pi r )\right)\,=\, \frac{\tilde{c}_{1,L,R}}{\sqrt{R+s}}\frac{\sqrt{2}}{\pi\cdot {\cos}(\pi R)}\,{\sin}(\pi s)
  \qquad\forall\,{R} \,\mbox{with\,}{\cos}(\pi R)\not=\,0\,.
\end{align*}
Now we show with sophisticated calculations that the limit of
${\lambda}_{1,S}({\cal A})$ 
as ${{\cal A}\to\infty}$ is ${\pi}^2$ as well.
\begin{corollary}
\label{PiSQuadratinf}
For the first eigenvalue ${\lambda}_{1,S}({\cal A})$ of the 
Stokes operator ${\boldsymbol S}$ (cf. Def. \ref{D5}) we observe that 
$\lim_{{\cal A}\to\infty} {\lambda}_{1,S}({\cal A})\,=\,{\pi}^2$.
\end{corollary}
\begin{proof}
We use  \eqref{trans2} and \eqref{Bess_DreiHalbe}.
For all
${\cal A}\,\in\,(0,\infty)$  the smallest positive solution of the
transcendental  equation \eqref{trans2} is called
$\kappa_{1,S}({\cal A})$
(cf. {also} \eqref{efSsigma_0}).
We abbreviate $\kappa \,=\,\kappa_{S}({\cal A})$
and $R\,=\,\frac{{\cal A}}{2}$. Then we get from \eqref{trans2} and
\eqref{Bess_DreiHalbe} that
\begin{align}
\label{transinf2}
0\,&=\,J_{\frac{3}{2}}(\kappa({1+R}))J_{-\frac{3}{2}}(\kappa {R})
\,-\,J_{\frac{3}{2}}(\kappa {R})J_{-\frac{3}{2}}(\kappa (1+{R})) \nonumber\\
\,&=\,   \frac{2}{\pi \kappa \sqrt{R(1+R)}}\cdot\,\big((-1)\left(-
{\cos(\kappa(1+{R}))}\,+\,
\frac{\sin(\kappa(1+{R}))}{\kappa(1+{R})}
\right)
\left(
{\sin(\kappa R)}\,+\,
\frac{\cos(\kappa {R})}{\kappa {R}}
\right)\,+\,
 \\
\,&\hspace*{2cm}+\,\left(
{\sin(\kappa(1+{R}))}\,+\,
\frac{\cos(\kappa(1+{R}))}{\kappa(1+{R})}
\right)
\left(
{-\cos(\kappa R)}\,+\,
\frac{\sin(\kappa {R})}{\kappa {R}}
\right)\big)\,.
{~} \nonumber
\end{align}
There are two ways to proceed: One can do the same calculations
with the Besselfunctions $J_{\frac{1}{2}}(.)$ and
$J_{-\frac{1}{2}}(.)$ (cf. \eqref{Bess_DreiHalbe}) and follow the proof
of Corollary \ref{PiQuadratinf}
or calculate \eqref{transinf2} directly.  We choose
the second one and conclude from \eqref{transinf2} that
\begin{align}
\label{calcinf2}
0\,=\,  \frac{2}{\pi \kappa \sqrt{R(1+R)}}\cdot\,\left(-\Big(1+\frac{1}{{\kappa}^2(1+{R})R}\Big)\sin(\kappa)\,+
\frac{1}{{\kappa}(1+{R})R} \,{\cos(\kappa)}
\right)\,.
\end{align}
Now we take into account the continuity of \eqref{transinf2} in
$\kappa$. By systematic asymptotic analysis we arrive
at ${\cal A}\,\to\,{\infty}$ or $R\,\to\,{\infty}$
like for the Nicholson's formula \cite[~4.1.3]{AAR} the (asymptotic) limiting value as the equation
\begin{align}
 0\,=\,  \,\frac{2}{\pi \kappa}
\sin(\kappa) \,\,,
\end{align}
with the smallest positive solution $\kappa\,=\,{\pi}$. This means $\lim_{{\cal A}\to\infty} {\lambda}_{1,S}({\cal A})\,=\,
\lim_{{\cal A}\to\infty} \kappa^2_{1,S}({\cal A})
\,=\,{\pi}^2$.
\end{proof}
\noindent Finally, we investigate the behaviour of the Poincar\'e
constant
as  ${{\cal A}\,\to\,\infty}$.
As for $c_{p}({\cal A})$ we see analogously,
that the (formal) limiting behaviour of $c_{p,S}({\cal A})$ is
\begin{align*}
\lim_{{\cal A}\,\to\,\infty}{\;\!} c_{p,S}({\cal A})={\frac{1}{\pi}}\,\,.
\end{align*}
Additionally we can specify an `approximate' limit
for the first
Stokes eigenfunctions with \eqref{Coordinate} for large ${\cal A}$
\begin{align}
{\underline{w}}_{1,{{\alpha}},S,{\cal A}}({\underline{x}})\,\approx\,{\tilde{c}_{{\;\!}\alpha}}\left(
 \,
{\sin}(\pi s) \cdot {\underline{\mathfrak{w}}}_{{\;\!}\alpha}(\vartheta,\varphi)\right)\,,\label{AppLimitSTEF}
\end{align}
where we have used the ${\underline{\mathfrak{w}}}_{{\;\!}\alpha}(\vartheta,\varphi)\,$ for $\alpha\,=\,-1,\,0,\,1$ in \eqref{eifstokes}.
The multiplicity of the eigenvalue $\pi^2$ of course remains constant
also for ${{\cal A}\,\to\,\infty}$. \\
In order to show \eqref{AppLimitSTEF} we use the $r-$dependent part of
\eqref{eifstokes} , $r\,:=\,\frac{\cal A}{2}+s\,=\,R+s \,$
\,{for $s\,\in\,[0,1]$}  and the abbreviation $\kappa_{1,S}({\cal A})\,:=\,\kappa $:
\begin{align}
f(r)\,=\,f(R+s)\,=:\,f(s)&=\frac{\tilde{C}_{\cal A}}{\sqrt{R+s}}
\left(
J_{\frac{3}{2}}({\kappa}(R+s))\,-\,
\frac{J_{{\frac{3}{2}}}({{\kappa R}})}
{J_{-{\frac{3}{2}}}({{\kappa R}})}
J_{-\frac{3}{2}}( \kappa (R+s)
)
\right)
\,.\label{eifstoAgross}
 \end{align}
We get by sophisticated calculations:
\begin{align}
f(s)&= \frac{\tilde{C}_{\cal A}}{R+s}
\sqrt{\frac{2}{\pi }}
\frac{\kappa R}{(\kappa  R\sin(\kappa R)+\cos(\kappa R))}
\cdot\,\left(\sin(\kappa s)\,+
\frac{\sin(\kappa s)-\kappa s \cos(\kappa s)}
{{\kappa}^2 (R+s) R} \,
\right)\,\,,
 \end{align}
where the relation \eqref{AppLimitSTEF} is received as asymptotic solution with $R\,\to\,\infty$ and Corollary
\ref{PiSQuadratinf} : $\lim_{{\cal A}\to\infty} {\kappa}_{1,S}({\cal
  A})\,=\,{\pi}\,=\,{\kappa}$.
\section{Results and Calculations}\label{Sec4}
For the numerical calculation of the Poincar\'e constants $c_{p,S}({\cal A})$ 
one has to take into account, that only their inverse values 
$\kappa_{1,S}({\cal A})$ are computable.
They are defined as roots of the  
transcendental equation
(\ref{trans2}). The use of ${\cal A}$ was proven and tested as an excellent scaling parameter
for the calculation of $\kappa_{1,S}({\cal A})$ 
over the range of ${\cal A}\,\in\,(0,\infty)$.  
The big advantage is a very small range for the roots of the
transcendental equation 
(\ref{trans2}).

\begin{figure}[th]
\centering
\scalebox{1.2}[1.6]{\includegraphics[width=0.55\textwidth]{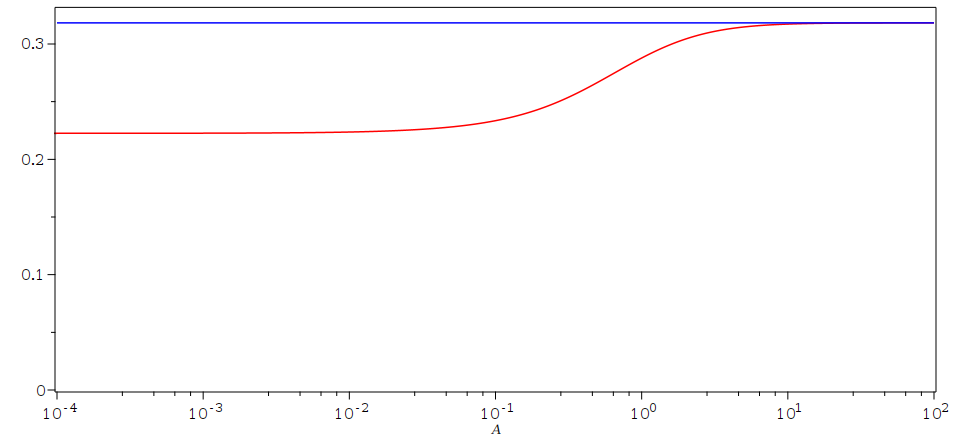}}
\caption{The constant value {\blue{$c_{p}({\cal
        A})\,=\,\frac{1}{\pi}$}} and calculated Stokes Poincar\'{e}
  constants 
{\red{$c_{p,S}({\cal A})$}} 
as functions of $\mathcal A$}
\label{fig:Stokes_plot}
\centering
  \includegraphics[width=0.73\textwidth]{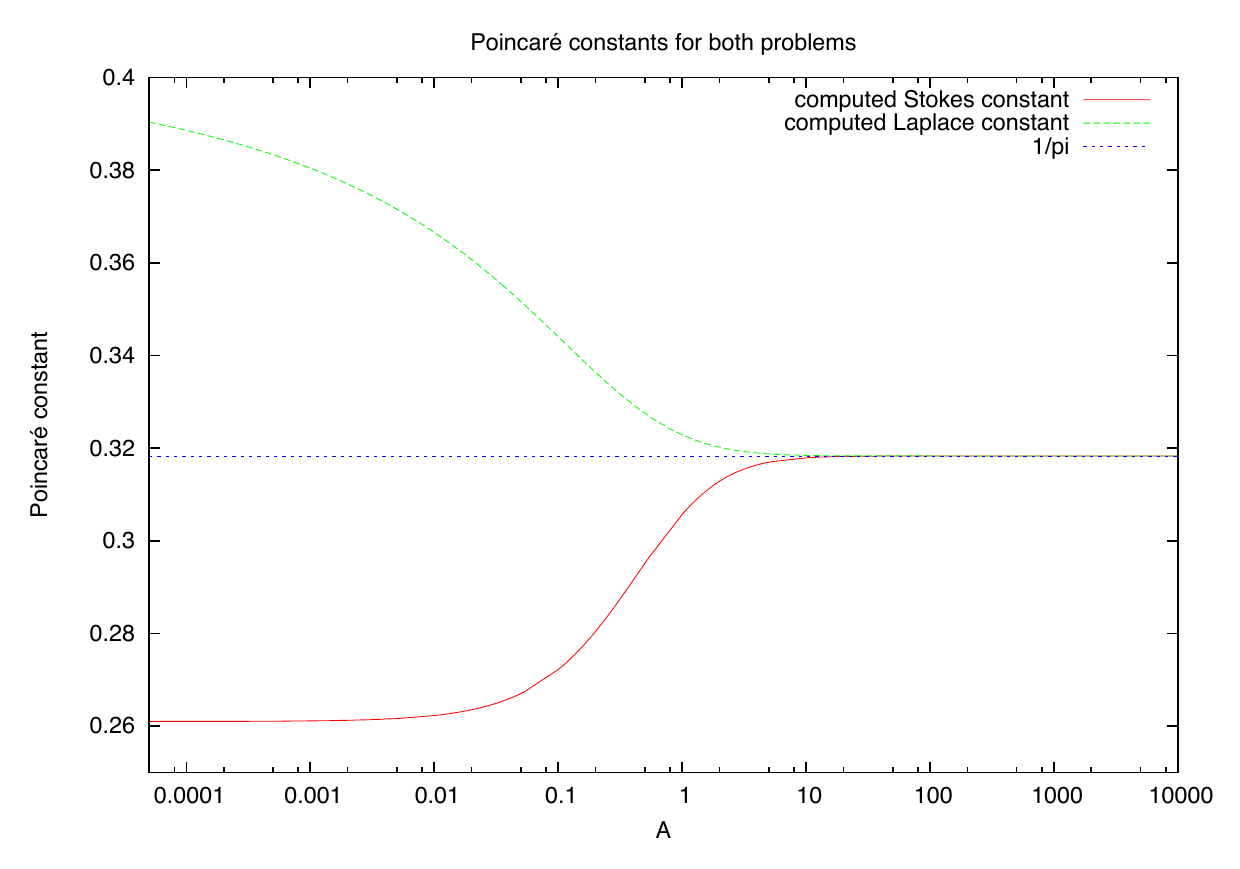}
\caption{Calculated Laplace and Stokes Poincar\'{e} constants in two
  dimensions 
as functions of ${\cal A}$}
\label{fig:Laplace_Po_plot}
\end{figure}

Standard Maple-procedures for the solution of transcendental equations
are used to find  
the first positive roots of (\ref{trans2})
and the corresponding Poincar\'e constants.
The first results of the calculations were written in data files as
${\cal A}$-pointwise exact 
$c_{p,S}({\cal A})$-values. The ${\cal A}$-grid
was chosen very fine
for small ${\cal A}$-values and rougher for large ${\cal A}$
because of the completely different behaviour
near ${\cal A}=0$ and  ${{\cal A}\,\to\,\infty}$, respectively. 
The calculations of {\red{$c_{p,S}({\cal A})$}} arre executed
with pointwise Maple 22 procedures. 
We use a logarithmic 
scale of ${\cal A}$ in the illustrations in Fig.
~\ref{fig:Stokes_plot}.
The graphical representation of the obtained result for Figure ~\ref{fig:Stokes_plot} was created with Maple.
The value of $c_{p,S}({\cal A})$ at ${\cal A}\,=\,0$ is 
${\red{c_{p,S}({0})\, = \,0.2225481584  }}$ .\\ [.2cm]

\noindent We also use Maple procedures for the calculations of
{\green{$c_{p}({\cal A})$}} and {\red{$c_{p,S}({\cal A})$}} in
Fig.~\ref{fig:Laplace_Po_plot} (cf. \cite{RuRuTh2016}). 
To facilitate the evaluation of the calculated results for the Poincar\'{e}  
constants we use a logarithmic 
scale for ${\cal A}$ in Fig.~\ref{fig:Laplace_Po_plot}.\pagebreak

\noindent Finally, let us state some results and some hypotheses for
the Laplace and 
    Stokes Poincar\'{e} constants (i.e. the first Laplace and Stokes
    eigenvalues, respectively)
    for the small-gap limit. 
    The limits of
{$c_{p}({\cal A})$} and {$c_{p,S}({\cal A})$} 
for ${\cal A} \to \infty$ are equal and 
\begin{align*}
\lim_{{\cal A}\,\to\,\infty}{\;\!} c_{p}({\cal A}) \,=\,\lim_{{\cal A}\,\to\,\infty}{\;\!} c_{p,S}({\cal A})={\frac{1}{\pi}}\,\,.
\end{align*}
This seems to be independent of the dimension of the annuli. 
This means
that this might be true for all $n$-dimensional annuli, $ 2\leq n <
\infty$. 
Moreover, for the first
eigenvalues $\lambda_{1,L}({\cal A})$ and 
$\lambda_{1,L}({\cal A})$ for any ${\cal A} \in (0,\infty)$. We only
have to check, that 
\begin{align*}
\lim_{{\cal A}\,\to\,\infty}{\;\!} \lambda_{1,L}({\cal A}) \,=\,{\pi}^2 \quad \text{or}\quad
\lim_{{\cal A}\,\to\,\infty}{\;\!} \lambda_{1,S}({\cal A}) \,=\,{\pi}^2 
\,\,.
\end{align*}
The multiplicity of the first eigenvalues $\lambda_{1,L}({\cal A})$
and $\lambda_{1,L}({\cal A})$ is constant and a continuous function of
${\cal A} \in (0,\infty)$. 
Finally, 
${\lambda}_{1,L}({\cal A})\,\leq\,{\lambda}_{1,S}({\cal A})$
independenly of the dimension of the annuli. Thus, $c_{p,S}({\cal A})\,\leq
\,c_{p}({\cal  A}) $ for all  ${\cal A} \in (0,\infty)$. 

\section*{Appendix}\label{App_Bes}
For the spherical Besselfunctions $J_{\frac{1}{2}}(t)$ and $J_{-\frac{1}{2}}(t)$
we can use 
(cf. \cite{Triebel}, 5.5.1, \cite{CouHil} or  \cite{Lewin})
\begin{eqnarray} \label{Bess_Halbe}
J_{\frac{1}{2}}(t)\,=\,\sqrt{\frac{2}{t\pi}}\cdot\,{\sin(t)} \quad  \quad ,\quad
J_{-\frac{1}{2}}(t)\,=\,\sqrt{\frac{2}{t\pi}}\cdot\,{\cos(t)} \,.
\end{eqnarray}
For 
$J_{\frac{3}{2}}(t)$ and $J_{-\frac{3}{2}}(t)$
the representation formulas by trigometric functions and spherical
Besselfunctions (cf.\cite{Lewin}) read
\begin{eqnarray} \label{Bess_DreiHalbe}
J_{\frac{3}{2}}(t)\,=\,\sqrt{\frac{2}{t\pi}}\cdot\,\left(-
{\cos(t)}\,+\,
\frac{\sin(t)}{t} 
\right)
\quad \,=\,
\frac{J_{\frac{1}{2}}(t)}{t} \,-\,J_{-\frac{1}{2}}(t)
\quad\nonumber \\
{~}\\
J_{-\frac{3}{2}}(t)\,=\,-\sqrt{\frac{2}{t\pi}}\cdot\,\left(
{\sin(t)}\,+\,
\frac{\cos(t)}{t} 
\right)
\quad \,=\,-\left(
\frac{J_{-\frac{1}{2}}(t)}{t} \,+\,J_{\frac{1}{2}}(t)\right)\,.
 \nonumber 
\end{eqnarray}
The Laplacian in the spherical coordinate system reads
\begin{align*} 
 \Delta_{r,\vartheta,\varphi} \, = \, \Big( 
\frac{\partial^{2}{ }}{\partial r^{2}} + \frac{2}{r}\frac{\partial { }}{\partial r} +\frac{1}{r^{2}}(
\frac{1}{\sin( \vartheta)}\frac{\partial{ }}{\partial \vartheta}{\sin( \vartheta)}\frac{\partial{ }}{\partial \vartheta}+
\frac{1}{\sin^{2}( \vartheta)}\frac{\partial^{2}{ }}{\partial \varphi^{2}})\Big)\,.
\end{align*}

\begin{thebibliography}{99}
\bibitem{ADN} S.~Agmon, A.~Douglis, and L.~Nirenberg,
Estimates near the boundary for solutions of elliptic 
partial differential equations satisfying general boundary conditions II,
Comm. Pure Appl. Math. {\bf 17}, 35-92  (1964).
\bibitem{AGBu}
T.~Akinaga , S.C.~Generalis , F.H.~Busse;
Tertiary and Quaternary States in the Taylor-Couette System
in Chaos, Solitons and Fractals
Volume 109, April 2018, Pages 107-117
https://doi.org/10.1016/j.chaos.2018.01.033
\bibitem{AAR} G.E. ~Andrews, R.~Askey, and R.~Roy,
Special Functions,
(Cambridge Univ.Press, Cambridge, New York, 1999).
\bibitem{Catbri} L.~Cattabriga,
{{Su un problema al contorno relativo si 
sistema di equazione di Stokes}},  
Rend. Mat. Univ. Padova {\bf 31},. 308-340 (1961).
\bibitem{CoFoi} P.~Constantin and C.~Foias,
{{Navier-Stokes Equations}}, 
(Univ.of Chic.Press, Chicago, 1988).	
\bibitem{CouHil} R.~Courant and  D.~Hilbert, 
{{Methoden der Mathematischen Physik}}, Vol.I and Vol.II
3. Aufl. (Springer, Berlin, Heidelberg, New York, 1968).
\bibitem{galdi1998}
{G.P.~Galdi}
An introduction to the mathematical theory of the Navier-Stokes
  equations, Vol.~1: Linearised steady problems
(Springer, New York, 1998).
\bibitem{GilTru} D. Gilbarg, N.S. Trudinger, Elliptic Partial Differential Equations of Second
Order, Grundlehren der mathematischen Wissenschaften 224, Reprint of the 1998
ed., (Springer, Berlin, Heidelberg, New York, 2001).
\bibitem{GirRav} V.~Girault and P.-A.~Raviart, 
{{Finite element approximation of the Navier-Stokes equations}}, 
(Springer, Berlin, 1979).
\bibitem{Jos}  D.D.~Joseph, 
{{Stability of Fluid Motions}},
Vol.I, (Springer, Berlin, Heidelberg, New York, 1976).
\bibitem{Junk} M.~Junk
{{Numerische Untersuchung der Stabilität der
Str\"omung im weiten Kugelspalt}}, (Cuvillier Verlag, Göttingen, 2005).  
\bibitem{KaiWahl}
R.~Kaiser, W.~von~Wahl, {{A New Functional for the Taylor-Couette Problem in the Small-Gap Limit}},
in Mathematical theory in fluid mechanics, Pitman Research Notes in Mathematics, Series 354, 
editors: G.P. Galdi, J. Malek, J. Necas, 114-134 (1996).
\bibitem{Lewin} W.~I.~Lewin und J.~I.~Grosberg
{Differentialgleichungen der mathematischen Physik}
(Verlag Technik, Berlin, 1952).
\bibitem{nazarov2000} A.I.~Nazarov,
The one-dimensional character of an extremum point of the
  {F}riedrichs inequality in spherical and plane layers,
Journal of Mathematical Sciences {\bf 102}, 5 (2000), 4473-4486.
\bibitem{passerini2009}
A.~Passerini, M.~R\r{u}\v{z}i\v{c}ka, and G.~Th{\"a}ter,
Natural convection between two horizontal coaxial cylinders,
ZAMM {\bf 89}, 5 (2009) 399-413.
\bibitem{passerini2010}
A.~Passerini, C.~Ferrario, M.~R\r{u}\v{z}i\v{c}ka, and G.~Th{\"a}ter,
Theoretical results on steady convective flows between horizontal
  coaxial cylinders,
SIAM Journal on Applied Mathematics {\bf 71}, 2 (2011) 465-486.
\bibitem{passerini2024}
A.~Passerini, B.~Rummler, M.~R\r{u}\v{z}i\v{c}ka, and G.~Th{\"a}ter,
Natural Convection in the Horizontal Annulus: Critical Rayleigh Number for the steady Problem,
ZAMM : Volume 105, Issue 3, March 2025, 
https://doi.org/10.1002/zamm.202300535
\bibitem{RumHab} B.~Rummler, 
{{Zur L\"osung der instation\"aren inkompressiblen 
Navier-Stokesschen Gleichungen in speziellen Gebieten}}, 
Magdeburg: Habilitation (1999/2000).
\bibitem{RumKug1} B.~Rummler, {The Eigenfunctions of the Stokes Operator in 
the open Unit Ball and in the open spherical Annulus},
Proc. of the  8th. Asian Computational Fluid Dynamics Conference, 
Hong Kong, 10-14 January, 2010
\bibitem{RumTh2024} B.~Rummler and G.~Th{\"a}ter, 
{{The Stokes Eigenvalue Problem on balls and annuli in three dimensions: Solutions with Poloidal
and Toroidal Fields}}, https://doi.org/10.48550/arXiv.2408.06948 (2024) 1-18
\bibitem{RuRuTh2016}
 B.~Rummler,, M.~R\r{u}\v{z}i\v{c}ka, and G.~Th{\"a}ter,
{Exact Poincar\'e constants in two-dimensional annuli}, ZAMM {\bf 97}, 1 (2017)
110–122.
\bibitem{Temam} R.~Temam, {{Navier-Stokes equations, theory and
      numerical 
analysis}}, 3rd edit., (North Holland, Amsterdam, 1984).
\bibitem{Triebel} H.~Triebel, {{Higher Analysis}}, 
  (Barth, Leipzig Berlin Heidelberg Amsterdam:, 1992).
\bibitem{Weid}
  J.~Weidmann, Stetige Abhängigkeit der Eigenwerte und
  Eigenfunktionen elliptischer Differentialoperatoren vom
  Gebiet. MATHEMATICA SCANDINAVICA, 54 (1984)
  51–69.
\end{thebibliography}
\end{document}